\input amstex\documentstyle{amsppt}  
\pagewidth{12.5cm}\pageheight{19cm}\magnification\magstep1
\topmatter
\title Nonsplit Hecke algebras and perverse sheaves\endtitle
\author G. Lusztig\endauthor
\address{Department of Mathematics, M.I.T., Cambridge, MA 02139}\endaddress
\thanks{Supported in part by National Science Foundation grant DMS-1303060 and by a Simons Fellowship.}
\endthanks
\endtopmatter   
\document

\define\pos{\text{\rm pos}}

\define\Bpq{\Bumpeq}

\define\bco{\bar{\co}}
\define\bcb{\bar{\cb}}
\define\bck{\bar{\ck}}

\define\uCS{\un{CS}}

\define\da{\dagger}
\define\dsv{\dashv}

\define\si{\sim}

\define\sqc{\sqcup}

\define\qua{\quad}

\define\hL{\hat L}
\define\hG{\hat G}

\define\bKK{\bar{\KK}}
\define\lb{\linebreak}

\define\op{\oplus}
   
\define\part{\partial}
\define\emp{\emptyset}
\define\imp{\implies}
\define\ra{\rangle}
\define\n{\notin}
\define\iy{\infty}
\define\m{\mapsto}
\define\do{\dots}
\define\la{\langle}

\define\sub{\subset}    

\define\T{\times}
\define\ti{\tilde}
\define\nl{\newline}
\redefine\i{^{-1}}

\define\un{\underline}

\define\ot{\otimes}
\define\bbq{\bar{\QQ}_l}

\define\Ad{\text{\rm Ad}}

\define\supp{\text{\rm supp}}

\define\bst{\bigstar}

\define\a{\alpha}
\redefine\b{\beta}
\redefine\c{\chi}
\define\g{\gamma}
\redefine\d{\delta}
\define\e{\epsilon}

\define\io{\iota}
\redefine\o{\omega}
\define\p{\pi}
\define\ph{\phi}
\define\ps{\psi}
\define\r{\rho}
\define\s{\sigma}
\redefine\t{\tau}
\define\th{\theta}
\define\k{\kappa}
\redefine\l{\lambda}
\define\z{\zeta}
\define\x{\xi}

\define\vt{\vartheta}

\redefine\G{\Gamma}

\define\Om{\Omega}

\define\Ph{\Phi}

\define\boc{\bold c}

\define\kk{\bold k}

\define\FF{\bold F}

\define\KK{\bold K}

\define\NN{\bold N}

\define\QQ{\bold Q}

\define\ZZ{\bold Z}

\define\ca{\Cal A}
\define\cb{\Cal B}

\define\cd{\Cal D}

\define\ch{\Cal H}
\define\ci{\Cal I}

\define\ck{\Cal K}
\define\cl{\Cal L}

\define\co{\Cal O}
\define\cp{\Cal P}

\define\cs{\Cal S}
\define\ct{\Cal T}
\define\cu{\Cal U}
\define\cv{\Cal V}
\define\cw{\Cal W}
\define\cz{\Cal Z}

\define\fK{\frak K}

\define\tc{\ti c}

\define\tf{\ti f}

\define\tj{\ti j}

\define\tp{\ti p}

\define\ty{\ti y}
\define\tA{\ti A}
\define\tB{\ti B}
\define\tC{\ti C}

\define\tP{\ti P}

\define\tU{\ti U}
\define\tV{\ti V}

\define\tZ{\ti Z}

\define\tKK{\ti{\KK}}
\define\sha{\sharp}

\define\bg{\bar g}

\define\ucs{\un{cs}}

\define\BE{Be}
\define\BBD{BBD}
\define\BRA{Br}
\define\GOR{Go}
\define\HE{He}
\define\KL{KL1}
\define\KLL{KL2}
\define\COX{L1}
\define\CLA{L2}
\define\ORA{L3}
\define\CLASSP{L4}
\define\HEC{L5}
\define\PCS{L6}
\define\CDGVIII{L7}
\define\PCSIII{L8}
\define\URCC{L9}

\head Introduction\endhead
\subhead 0.1\endsubhead 
Let $\cw$ be a Coxeter group such that the set $K$ of simple reflections of $\cw$ is finite. Let 
$l:\cw@>>>\NN$ be the length function and let $\cl:\cw@>>>\NN$ be a weight function that is, a 
function such that $\cl(xy)=\cl(x)+\cl(y)$ whenever $x,y\in\cw$ satisfy $l(xy)=l(x)+l(y)$. Let 
$\ca=\ZZ[v,v\i]$ where $v$ is an indeterminate. The Hecke algebra of $\cw$ (relative to $\cl$) is the free 
$\ca$-module $\ch$ with basis $\{\ct_z;z\in W\}$ with the associative algebra structure defined by the rules
$(\ct_z+v^{-\cl(z)}(\ct_z-v^{\cl(z)})=0$ if $z\in K$ and $\ct_{xy}=\ct_x\ct_y$ whenever $x,y\in\cw$ 
satisfy $l(xy)=l(x)+l(y)$. Let $\{c_z;z\in\cw\}$ be the basis of $\ch$ defined in \cite{\HEC, 5.2} in terms 
of $\cw,\cl$. We have $c_z=\sum_{t\in\cw}p_{t,z}\ct_t$ where $p_{t,z}\in\ZZ[v\i]$ is zero for all but 
finitely many $t$. In the case where $\ch$ is {\it split} that is $\cl=l$, we have
$p_{t,z}(v)=v^{-l(z)+l(t)}P_{t,z}(v^2)$ where $P_{t,z}$ is the polynomial defined in \cite{\KL}; in this
case, $p_{t,z}$ can be interpreted geometrically in terms of intersection cohomology (at least in the
crystallographic case), see \cite{\KLL}, and this interpretation has many interesting consequences.

In this paper we are interested in the case where $\ch$ is not split (and not even quasisplit, in the sense 
of \cite{\HEC, 16.3}).
 As shown in \cite{\COX}, \cite{\CLA} (resp. in \cite{\CLASSP}), such $\ch$, with $\cw$
a finite (resp. affine) Weyl group, appear as endomorphism algebras of representations of Chevalley groups 
over $\FF_q$ (resp. $\FF_q((\e))$) induced from unipotent cuspidal representations of a Levi quotient of a 
parabolic (resp. parahoric) subgroup. Our main goal is to describe the elements $p_{t,z}$ coming from 
(nonsplit) $\ch$ which appear in this way in representation theory, in geometric terms, involving perverse 
sheaves. In this paper we outline a strategy to achieve this goal using geometry
based on the theory of parabolic character sheaves \cite{\PCS}.
For simplicity we focus on Chevalley groups over $\FF_q$, but a similar strategy works for Chevalley
groups over $\FF_q((\e))$ in which case the affine analogue of parabolic character sheaves \cite{\PCSIII}
should be used. We illustrate in detail our strategy in the special case where $\ch$ is the endomorphism 
algebra of 
the representation of $SO_9(\FF_q)$ induced from a unipotent cuspidal representation of a Levi quotient of 
type $B_2$ of a parabolic subgroup. (This is the smallest example in which some 
$p_{t,z}$ can be outside $\NN[v\i]$. This $\ch$ is of type $B_2$ with $\cl$ taking the values $1$ and $3$ at
the two simple reflections.) 
The main effort goes into computing as much as possible of the cohomology sheaves of parabolic
character sheaves in this case. (For this we use the complete knowledge of the polynomials $P_{y,w}$ for the
Weyl group of type $B_4$, the knowledge of the multiplicity formulas for unipotent character sheaves on 
$SO_5$ and some additional arguments.) Eventually the various $p_{t,z}$ can be reconstructed from the 
information contained in the various cohomology sheaves of parabolic character sheaves. I expect that 
similar results hold for all $\ch$ appearing as above from representation theory.
(A conjecture in this 
direction is formulated in 3.11. It makes more precise a conjecture that was stated in \cite{\HEC, 27.12} 
before the theory of parabolic character sheaves was available.) 

\subhead 0.2. Notation\endsubhead
If $X$ is a finite set, $\sha X$ denotes the number of elements of $X$.
Let $2^X$ be the set of subsets of $X$.
If $\G$ is a group, $\G'$ is a subset of $\G$ and $\g\in\G$ we set ${}^\g\G'=\g\G\g\i$.
Let $\kk$ be an algebraically closed field. All algebraic varieties are assumed to be over $\kk$.
If $X$ is an algebraic variety, $\cd(X)$ denotes the bounded derived category of $\bbq$-sheaves on $X$ 
($l$ is a fixed prime number invertible in $\kk$). We will largely follow the notation of \cite{\BBD}.
If $K\in\cd(X)$ and $A$ is a simple perverse sheaf on $X$ we write $A\dsv K$ instead of: ''$A$ is a
composition factor of $\op_{j\in\ZZ}{}^pH^j(K)$''.
For $K\in\cd(X)$, $n\in\ZZ$ we write $K\la n\ra$ instead of 
$K[n](n/2)$ (shift, folowed by Tate twist). 
If $f:X@>>>Y$ is a smooth morphism of algebraic varieties all of whose fibres are irreducible of dimension
$d$ and $K\in\cd(Y)$, we set $f^\bst(K)=f^*(K)\la d\ra$.  
If $X$ is an irreducible algebraic variety, $|X|$ denotes the dimension of $X$. 
For any connected affine algebraic group $H$, $U_H$ denotes the unipotent radical of $H$.
If $\kk$ is an algebraic closure of a finite field $\FF_q$ and $X$ is an algebraic variety over $\kk$ with a
fixed $\FF_q$-structure, we shall denote by $\cd_m(X)$ the bounded derived category of mixed 
$\bbq$-complexes on $X$.

Assume that $C\in\cd(X)$ and that $\{C_i;i\in\ci\}$ is a family of objects of $\cd(X)$. We shall write
$C\Bpq\{C_i;i\in\ci\}$ if the following condition is satisfied: there exist distinct elements
$i_1,i_2,\do,i_s$ in $\ci$, objects $C'_j\in\cd(X)$ ($j=0,1,\do,s$) and distinguished triangles 
$(C'_{j-1},C'_j,C_{i_j})$ for $j=1,2,\do,s$ such that $C'_0=0$, $C'_s=C$; moreover, $C_i=0$ unless $i=i_j$ 
for some $j\in[1,s]$. The same definition can be given with $\cd(X)$ replaced by $\cd_m(X)$.

\head 1. The variety $Z_J$ and its pieces ${}^wZ_J$\endhead
\subhead 1.1\endsubhead
We fix an affine algebraic group $\hG$ whose identity component $G$ is reductive; we also fix a connected 
component $G^1$ of $\hG$. Let $\cb$ be the variety of Borel subgroups of $G$. Recall that the Weyl group $W$
of $G$ naturally indexes the set of $G$-orbits of the simultaneous conjugation action of $G$ on $\cb\T\cb$. 
We write $\co_w$ for the $G$-orbit indexed by $w\in W$. Let $l:W@>>>\NN$ be the length function 
$w\m l(w)=|\co_w|-|\co_1|$; let $I=\{w\in W;l(w)=1\}$. Recall that $(W,I)$ is a finite Coxeter group. Let 
$\le$ be the standard partial order on $W$. For any $J\sub I$ let $W_J$ be the subgroup of $W$ generated by 
$J$; let $W^J$ (resp. ${}^JW$) be the set of all $w\in W$ such that $l(wa)=l(w)+l(a)$ (resp. 
$l(aw)=l(a)+l(w)$) for any $a\in W_J$. For $J\sub I$, $K\sub I$, let ${}^KW^J={}^KW\cap W^J$. Any 
$(W_K,W_J)$-double coset $X$ in $W$ contains a unique element of ${}^KW^J$ denoted by $\min(X)$.

We define an automorphism $\d:W@>>>W$ by 
$$(B,B')\in\co_w,g\in G^1\imp({}^gB,{}^gB')\in\co_{\d(w)}.$$
We have $\d(I)=I$. For $J\sub I$ we write ${}^\d J$ instead of $\d(J)$; let $N_{J,\d}=\{w\in W;wJw\i
={}^\d J\}$. We have $N_{J,\d}\sub{}^{{}^\d J}W^J$.   

\subhead 1.2\endsubhead
Let $\cp$ be the set of parabolic subgroups of $G$. For $P\in\cp$ we set $\cb_P=\{B\in\cb;B\sub P\}$. For 
any $J\sub I$ let $\cp_J$ be the set of all $P\in\cp$ such that 
$$\{w\in W;(B,B')\in\co_w\text{ for some }B,B'\in\cb_P\}=W_J.$$
If $J\sub I$ and $P\in\cp_J,g\in G^1$ then ${}^gP\in\cp_{{}^\d J}$. We have $\cp=\sqc_{J\sub I}\cp_J$, 
$\cp_\emp=\cb$, $\cp_I=\{G\}$. For $B\in\cb$, $J\sub I$, there is a unique $P\in\cp_J$ such that $B\sub P$; 
we set $P=B_J$. For $J,K\sub I$, $P\in\cp_K$, $Q\in\cp_J$, the set
$$\{w\in W;(B,B')\in\co_w\text{ for some }B\in\cb_P,B'\in\cb_Q\}$$
is a single $(W_K,W_J)$-double coset in $W$, hence it contains a unique element in ${}^KW^J$ denoted by 
$\pos(P,Q)$. We set $P^Q=(P\cap Q)U_P$. We have 
$P^Q\in\cp_{K\cap\Ad(u)J}$ where $u=\pos(P,Q)$ and $U_{P^Q}=U_P(P\cap U_Q)$ hence 
$$|U_{P^Q}|=|P\cap U_Q|+|U_P|-|U_P\cap U_Q|.$$
Note that the condition that $P,Q$ contain a common Levi subgroup is equivalent to the condition that
$Ad(u\i)K=J$; in this case we have $P^Q=P$, $Q^P=Q$.

\subhead 1.3\endsubhead
Let $J\sub I$. Following B\'edard \cite{\BE} (see also \cite{\PCS, 2.4, 2.5}), for any $w\in{}^{{}^\d J}W$ 
we define a sequence $J=J_0\supset J_1\supset J_2\supset\do$ in $2^I$ and a sequence $w_0,w_1,\do$ in $W$ 
inductively by

$J_0=J$, 

$J_n=J_{n-1}\cap\d\i(\Ad(w_{n-1}J_{n-1})$, for $n\ge1$,

$w_n=\min(W_{{}^\d J}wW_{J_n})$, for $n\ge0$.
\nl
For $n\gg0$ we have $J_n=J_{n+1}=\do$ and $w_n=w_{n+1}=\do$; we set $J_\iy^w=J_n$ for $n\gg0$ and 
$w_\iy=w_n$ for $n\gg0$. 

According to {\it loc.cit.}, this gives a bijection $\k:{}^{{}^\d J}W@>\si>>\ct(J,\d)$ where $\ct(J,\d)$ is 
the set consisting of all sequences $(J_n,w_n)_{n\ge0}$ in $2^I\T W$ where 

$J_n=J_{n-1}\cap\d\i(\Ad(w_{n-1})J_{n-1})$ for $n\ge1$,

$w_n\in{}^{{}^\d J_n}W^{J_n}$ for $n\ge0$,     

$w_n\in w_{n-1}W_{J_{n-1}}$ for $n\ge1$.
\nl
The inverse bijection $\ct(J,\d)@>>>{}^{{}^\d J}W$ is given by $(J_n,w_n)_{n\ge0}\m w_\iy$. 

Assume that $w\in{}^{{}^\d J}W$. We have $J^w_\iy=\d\i(\Ad(w)J^w_\iy)$. 
Hence there is a well defined Coxeter group automorphism $\t_w:W_{{}^\d(J^w_\iy)}@>>>W_{{}^\d(J^w_\iy)}$ 
given by $x\m\t_w(x)=w\d\i(x)w\i$.

\subhead 1.4\endsubhead
Let $J\sub I$. We set
$$\tZ_J=\{(P,P',g)\in\cp_J\T\cp_{{}^\d J}\T G^1;{}^gP=P'\},$$ 
$$Z_J=\{(P,P',gU_P);P\in\cp_J,P'\in\cp_{{}^\d J},gU_P\in G^1/U_P,{}^gP=P'\}.$$

The variety $Z_J$ is the main object of this paper. (In \cite{\PCS, 3.3} $\tZ_J,Z_J$ are denoted by 
$\cz_{J,\d}$, $Z_{J,\d}$.)

Now $G\T G$ acts on $\tZ_J$ by $(h,h'):(P,P',g)\m({}^hP,{}^{h'}P',h'gh\i)$ and on $Z_J$ by 
$(h,h'):(P,P',gU_P)\m({}^hP,{}^{h'}P',h'gh\i U_{{}^hP})$. In this paper we shall restrict these actions to 
$G$ viewed as the diagonal in $G\T G$. Let $e_J:\tZ_J@>>>Z_J$ be the obvious map (an affine space bundle).

Following \cite{\PCS} we will define a partition of $Z_J$ into pieces indexed by ${}^{{}^\d J}W$.

To any $(P,P',g)\in\tZ_J$ we associate an element $w_{P,P',g}\in W$ by the requirements (i),(ii) below. (We 
set $z=\pos(P',P)\in{}^{{}^\d J}W^J$.)

(i) $w_{P,P',g}=w_{P_1,P'_1,g}$ where $P_1=P^{{}^{g\i}P}=({}^{g\i}P\cap P)U_P\in\cp_{J\cap\d\i(\Ad(z)J)}$,
$P'_1=P'{}^P\in\cp_{{}^\d J\cap\Ad(z)J}$;

(ii) $w_{P,P',g}=z$ if $z\in N_{J,\d}$.
\nl
These conditions define uniquely $w_{P,P',g}$ by induction on $\sha J$. If $z\in N_{J,\d}$ (in particular if
$\sha J=0$) then $w_{P,P',g}$ is given by (ii) (and (i) is satisfied since $P_1=P$). If $z\n N_{J,\d}$,
then $\sha(J\cap\d\i(\Ad(z)J))<\sha J$ and $w_{P,P',g}$ is determined by (i) since $w_{P_1,P'_1,g}$ is known
from the induction hypothesis.

From the definitions we see that the map $\tZ_J@>>>W$, $(P,P',g)\m w_{P,P',g}$ is the composition
$\tZ_J@>>>\ct(J,\d)@>\k\i>>{}^{{}^\d J}W$ (the first map is as in \cite{\PCS, 3.11}); in particular for any 
$(P,P',g)\in\tZ_J$ we have $w_{P,P',g}\in{}^{{}^\d J}W$, $w_{P,P',g}\in W_{{}^\d J}\pos(P',P)W_J$. For any 
$w\in{}^{{}^\d J}W$, we set 
$$\align&{}^w\tZ_J=\{(P,P',g)\in\tZ_J;w_{P,P',g}=w\},\\&
{}^wZ_J=\{(P,P',gU_P)\in Z_J;w_{P,P',g}=w\}.\endalign$$ 
The subsets $\{{}^w\tZ_J;w\in{}^{{}^\d J}W\}$ are said to be the pieces of $\tZ_J$; they form a partition of
$\tZ_J$. The subsets 
$\{{}^wZ_J;w\in{}^{{}^\d J}W\}$ are said to be the pieces of $Z_J$; they form a partition of $Z_J$. We have 
${}^w\tZ_J=e_J\i({}^wZ_J)$, ${}^wZ_J=e_J({}^w\tZ_J)$,
and ${}^w\tZ_J,{}^wZ_J$ are stable under the $G$-actions on $\tZ_J,Z_J$.

\subhead 1.5\endsubhead
Let $z\in{}^{{}^\d J}W^J$ and let $\Om=W_{{}^\d J}zW_J$. We set $J_1=J\cap\d\i(\Ad(z)J)$. Let 
$$\tZ_{J,\Om}=\{(P,P',g)\in\tZ_J;\pos(P',P)=z\},$$
a locally closed subvariety of $\tZ_J$. Let 
$$\tZ^\da_{J_1,\Om}=\{(Q,Q',g)\in\tZ_{J_1};\pos(Q',Q)\in zW_J\},$$
a locally closed subvariety of $\tZ_{J_1}$. By \cite{\PCS, 3.2}, 
$$\ti\a:\tZ_{J,\Om}@>>>\tZ^\da_{J_1,\Om},\qua
(P,P',g)\m(P_1,P'_1,g)\text{ where }P'_1=P'{}^P,P_1={}^{g\i}P'_1\tag a$$
is a well defined morphism; by \cite{\PCS, 3.6}, 

(b) {\it $\ti\a$ is an isomorphism.}
\nl
Let 
$$Z_{J,\Om}=\{(P,P',gU_P)\in Z_J;\pos(P',P)=z\},$$
a locally closed subvariety of $Z_J$. Let 
$$Z^\da_{J_1,\Om}=\{(Q,Q',gU_Q)\in Z_{J_1};\pos(Q',Q)\in zW_J\},$$ 
a subvariety of $Z_{J_1}$. Now $\ti\a$ induces a morphism $\a:Z_{J,\Om}@>>>Z^\da_{J_1,\Om}$ given by
$$(P,P',gU_P)\m(P_1,P'_1,gU_{P'}\text{ where }P'_1=P'{}^P,P_1={}^{g\i}P'_1.\tag c$$
From \cite{\PCS, 3.7-3.10} we see that 

(d) {\it $\a$ is an affine space bundle with fibres of dimension $|U_P\cap P'|-|U_P\cap U_{P'}|$ for 
some/any $(P,P')\in\cp_J\T\cp_{{}^\d J}$ such that $\pos(P',P)=z$.}
\nl
Next we note that, for $w\in{}^{{}^\d J}W$ such that $w\in\Om$, we have 

${}^w\tZ_J\sub\tZ_{J,\Om}$, ${}^wZ_J\sub Z_{J,\Om}$, ${}^w\tZ_{J_1}\sub\tZ^\da_{J_1,\Om}$, 
${}^wZ_{J_1}\sub Z^\da_{J_1,\Om}$.
\nl
(We use that $w\in zW_J$.) Using the definitions we have

(e) ${}^w\tZ_J=\ti\a\i({}^w\tZ_{J_1}), {}^wZ_J=\a\i({}^wZ_{J_1})$.
\nl
Moreover, using (b),(d) we deduce:

(f) {\it $\ti\a$ restricts to a bijection $\ti\vt_{J,w}:{}^w\tZ_J@>>>{}^w\tZ_{J_1}$;}

(g) {\it $\a$ restricts to a map $\vt_{J,w}:{}^wZ_J@>>>{}^wZ_{J_1}$ all of whose fibres are affine spaces of 
dimension $|U_P\cap P'|-|U_P\cap U_{P'}|$ for some/any $(P,P')\in\cp_J\T\cp_{{}^\d J}$ such that 
$\pos(P',P)=z$.}

\proclaim{Proposition 1.6}Let $J\sub I$, $w\in{}^{{}^\d J}W$. Define $z\in{}^{{}^\d J}W^J$ by $w\in zW^J$.

(a) ${}^w\tZ_J$ (resp. ${}^wZ_J$) is a smooth irreducible locally closed subvariety of $\tZ_J$ (resp. $Z_J$)
of dimension $l(w)+|G|$ (resp. $l(w)+|G|-|\cp_J|$).

(b) Let $J_1=J\cap\d\i(\Ad(z)J)$. Then $\ti\vt_{J,w}:{}^w\tZ_J@>>>{}^w\tZ_{J_1}$ (see 1.5(f)) is an 
isomorphism and $\vt_{J,w}:{}^wZ_J@>>>{}^wZ_{J_1}$ is a smooth morphism all of whose fibres are affine 
spaces of dimension $|\cp_{J_1}|-|\cp_J|$.
\endproclaim
Assume first that $z\in N_{J,\d}$. Then $z=w$ and 
$${}^w\tZ_J=\{(P,P',g)\in\tZ_J;\pos(P',P)=z\}$$
(resp. ${}^wZ_J=\{(P,P',gU_P)\in Z_J;\pos(P',P)=z\}$) is the inverse image of 
$${}^zY_J=\{(P,P')\in\cp_J\T\cp_{{}^\d J};\pos(P',P)=z\}$$
under the obvious map $\tZ_J@>>>\cp_J\T\cp_{{}^\d J}$ (resp. $Z_J@>>>\cp_J\T\cp_{{}^\d J}$), a smooth map 
with fibres isomorphic to $P$ (resp. $P/U_P$) for $P\in\cp_J$. Since ${}^zY_J$ is smooth, irreducible, 
locally closed in $\cp_J\T\cp_{{}^\d J}$, of dimension $l(z)+|\cp_J|$, it follows that in this case 
${}^w\tZ_J$ (resp. ${}^wZ_J$) is a smooth, irreducible, locally closed subvariety of $\tZ_J$ (resp. $Z_J$) 
of dimension $l(z)+|\cp_J|+|P|$ (resp. $l(z)+|\cp_J|+|P/U_P|$). Thus (a) follows in this case.

If $z\n N_{J,\d}$ then, setting $J_1=J\cap\d\i(\Ad(z)J)$, we have $\sha J_1<\sha J$ and we can assume that
(a) holds when $J$ is replaced by $J_1$. Using 1.5(e),(b),(d) we deduce that ${}^w\tZ_J$ (resp. ${}^wZ_J$) 
is a smooth irreducible locally closed subvariety of $\tZ_J$ (resp. $Z_J$) of dimension $l(w)+|G|$ (resp. 
$l(w)+|G|-|\cp_{J_1}|+|U_P\cap P'|-|U_P\cap U_{P'}|$ for some/any $(P,P')\in{}^zY_J$). To complete the proof
of (a) it is then enough to note that $|{}^wZ_J|=|{}^w\tZ_J|-|\cp_J|$. The proof above shows that 
$$|\cp_{J_1}|-|\cp_J|=|U_P\cap P'|-|U_P\cap U_{P'}|$$
for some/any $(P,P')\in{}^zY_J$.

Now (b) follows immediately from (a) using 1.5(e),(b),(d).

\subhead 1.7. Examples\endsubhead
In the case where $J=I$ we can identify $Z_J=G^1$. It has a unique piece, ${}^1Z_J=Z_J$. In the case where
$J=\emp$, $Z_J$ is a torus bundle over $\cb^2$. The pieces of $Z_J$ are the inverses images of the
$G$-orbits $\co_w$ ($w\in W$) under $Z_J@>>>\cb^2$.

Assume now that $\hG=G=GL(V)$ where $V$ is a $\kk$-vector space of dimension $3$.
We can identify the projective space $P(V)$ with $\cp_J$ for a certain $1$-element subset $J$ of $I$.
Then $\tZ_J$ becomes the set of triples $(L,L',g)$ where $L,L'$ are lines in $V$ and 
$g\in GL(V)$ carries $L$ to $L'$; since $L'$ is determined by $L,g$ we can identify
$\tZ_J$ with the set of pairs $(L,g)$ where $L$ is a line in $V$ and $g\in GL(V)$.
For any $r\in[1,3]$ let ${}^r\tZ_J$ be the set of all $(L,g)\in\tZ_J$ such that $L,gL,g^2L,\do$
span an $r$-dimensional subspace. 
Then the subsets ${}^r\tZ_J$ ($r\in[1,3]$) are exactly the pieces of $\tZ_J$.
Now $Z_J$ becomes the set of quadruples $(L,L',\g,\ti\g)$ where $L,L'$ are lines in $V$ and 
$\g:L@>\si>>L'$, $\ti\g:V/L@>\si>>V/L'$ are isomorphisms of vector spaces.
Let ${}^1Z_J=\{(L,L',\g,\ti\g)\in Z_J;L=L'\}$.
Let ${}^2Z_J$ be the set of all $(L,L',\g,\ti\g)\in Z_J$ such that $L\ne L'$
and $\ti\g$ carries the image of $L'$ in $V/L$ to the image of $L$ in $V/L'$.
Let ${}^3Z_J$ be the set of all $(L,L',\g,\ti\g)\in Z_J$ such that $L\ne L'$ and such that the following 
holds: denoting by $L_1$ the image of $L'$ in $V/L$, by $L'_1$ the image of $L$ in $V/L'$, 
and setting $L_2=\ti\g\i(L'_1)\sub V/L$, $L'_2=\ti\g(L_1)\sub V/L'$, we have $L_2\ne L_1$ (hence
$L'_2\ne L'_1$) and $\ti\g(L_1\op L_2)=L'_1\op L'_2$.
Then ${}^1Z_J,{}^2Z_J,{}^3Z_J$ are exactly the pieces of $Z_J$.

Assume now that $\hG=G=Sp(V)$ where $V$ is a $\kk$-vector space of dimension $4$ with a given 
nondegenerate symplectic form $\la,\ra$. We can identify the projective space $P(V)$ with $\cp_J$ for a 
certain $1$-element subset $J$ of $I$.
Then $\tZ_J$ becomes the set of triples $(L,L',g)$ where $L,L'$ are lines in $V$ and $g\in Sp(V)$
carries $L$ to $L'$; since $L'$ is determined by $L,g$ we can identify
$\tZ_J$ with the set of pairs $(L,g)$ where $L$ is a line in $V$ and $g\in Sp(V)$.

For $r\in[1,2]$ let ${}^r\tZ_J$ be the set of all $(L,g)\in\tZ_J$ such that $L,gL,g^2L,\do$
span an $r$-dimensional isotropic subspace of $V$ and let ${}^rZ_J'$ be 
the set of all $(L,g)\in\tZ_J$ such that 
$\la L,gL\ra=\la L,g^2L\ra=\do=\la L,g^{r-1}L\ra=0$, $\la L,g^rL\ra\ne0$.
Then the subsets ${}^r\tZ_J,{}^r\tZ'_J$ ($r=1,2$) are exactly the pieces of $\tZ_J$.

Now $Z_J$ becomes the set of quadruples $(L,L',\g,\ti\g)$ where $L,L'$ are lines in $V$ and $\g:L@>\si>>L'$ 
is an isomorphism of vector spaces and $\ti\g:L^\perp/L@>\si>>L'{}^\perp/L'$ is a symplectic isomorphism. 
Let ${}^1Z_J=\{(L,L',\g,\ti\g)\in Z_J;L=L'\}$. Let ${}^2Z_J$ be the set of all $(L,L',\g,\ti\g)\in Z_J$ such
that $L\ne L',\la L,L'\ra=0$ and $\ti\g$ carries the image of $L'$ in $L^\perp/L$ to the image of $L$ in 
$L'{}^\perp/L'$. Let ${}^2Z'_J$ be the set of all $(L,L',\g,\ti\g)\in Z_J$ such that $L\ne L'$, 
$\la L,L'\ra=0$ 
and $\ti\g$ does not carry the image of $L'$ in $L^\perp/L$ to the image of $L$ in $L'{}^\perp/L'$.
Let ${}^1Z'_J=\{(L,L',\g,\ti\g)\in Z_J;\la L,L'\ra\ne0\}$.
Then the subsets ${}^1Z_J,{}^2Z_J,{}^2Z'_J,{}^1Z'_J$ are exactly the pieces of $Z_J$.

\subhead 1.8\endsubhead
Let $J\sub I$. Let 
$$\cb^2_J=\{(B,B',gU_{B_J});(B,B')\in\cb^2,gU_{B_J}\in G^1/U_{B_J},{}^gB=B'\};$$ 
this is well defined since $U_{B_J}\sub B$. Define $\p_J:\cb^2_J@>>>Z_J$ by 
$$(B,B',gU_{B_J})\m(B_J,B'_{{}^\d J},gU_{B_J}).$$

For any $y\in W$ let 
$$\cb^2_{J,y}=\{(B,B',gU_{B_J})\in\cb^2_J;(B,B')\in\co_y\}$$
and let $\p_{J,y}:\cb^2_{J,y}@>>>Z_J$ be the restriction of $\p_J$. 
The statements (a),(b) below can be deduced from \cite{\PCS, 4.14}.

(a) {\it For any $w\in{}^{{}^\d J}W$, the image of $\p_{J,w\i}:\cb^2_{J,w\i}@>>>Z_J$ is exactly ${}^wZ_J$.}

(b) {\it If $w\in{}^{{}^\d J}W$ and $x\in W_{{}^\d(J^w_\iy)}$ then the image of 
$\p_{J,w\i x}:\cb^2_{J,w\i x}@>>>Z_J$ is contained in ${}^wZ_J$.}
\nl
Note that (a) gives an alternative description of ${}^wZ_J$ as the image of $\p_{J,w\i}$.

\subhead 1.9\endsubhead
Let $J\sub I$ and let $w\in{}^{{}^\d J}W$. In \cite{\HE, 4.6} it is shown that the closure of ${}^wZ_J$ in
$Z_J$ is equal to $\cup_{w'\in{}^{d(J)}W;w'\le_Jw}{}^{w'}Z_J$ where $w'\le_Jw$ is defined by the condition:
$\d(u)w'u\i\le w$ for some $u\in W_J$.

\head 2. Unipotent character sheaves on $Z_J$ and on ${}^wZ_J$\endhead
\subhead 2.1\endsubhead
Let $J\sub I$. For $y\in W$ we set $\KK^y_J=(\p_{J,y})_!\bbq\in\cd(Z_J)$. A {\it unipotent character sheaf}
on $Z_J$ is by definition a simple perverse sheaf $A$ on $Z_J$ such that $A\dsv\KK^y_J$ for some $y\in W$. 
This is a special case of what in \cite{\PCS} is referred to as a {\it parabolic character sheaf}.
Let $CS(Z_J)$ be the collection of unipotent character sheaves on $Z_J$.

In the case where $J=I$, we can identify $Z_J=Z_I$ with $G^1$. Hence there is a well defined notion of 
unipotent character sheaf on $G^1$. In this case, for $y\in W$ we have
$$\cb^2_{I,y}=\{(B,B',g)\in\cb\T\cb\T G^1;B'={}^gB,(B,B')\in\co_y\}$$
and $\p_{I,y}:\cb^2_{J,y}@>>>Z_I$ is given by $(B,B',g)\m g$. 

In the case where $J=\emp$, for any $y\in W$, $\p_{J,y}$ is the inclusion of 
$$\{(B,B',gU_B);(B,B')\in\co_y,g\in G^1,{}^gB=B'\}$$ into 
$$Z_\emp=\{(B,B',gU_B);(B,B')\in\cb^2,g\in G^1,{}^gB=B'\}.$$ 
It follows that $CS(Z_\emp)$ consists of the simple 
perverse sheaves on $Z_\emp$ which are (up to shift) the inverse images under $Z_\emp@>>>\cb^2$ of 
the simple $G$-equivariant perverse sheaves on $\cb^2$.

\subhead 2.2\endsubhead
Let $J\sub I$ and let $w\in{}^{{}^\d J}W$. We define a collection of simple perverse sheaves $CS({}^wZ_J)$ 
on ${}^wZ_J$ (said to be unipotent character sheaves on ${}^wZ_J$) by induction on $\sha J$ as follows. We 
set $z=\min(W_{{}^\d J}wW_J)$. 

Assume first that $z\in N_{J,\d}$ so that $z=w$.
For any $(P,P')\in{}^zY_J$ we denote by $\cs_{P,P'}$ the set of common Levi subgroups of $P,P'$; this is a
nonempty set. For any $L\in\cs_{P,P'}$ we denote by $\hL$ the normalizer of $L$ in $\hG$.
Note that the identity component of $\hL$ is $L$. We set $\hL^1=\{g\in G^1;{}^gP=P',{}^gL=L\}\sub\hL$. If 
$g,g'\in\hL^1$ then, setting $g'=gh$, we have $h\in P\cap\hL$ hence $h\in L$; we see that $\hL^1$ is a single
connected component of $\hL$. We have a diagram
$$\hL^1@<c<<G\T\hL^1@>c'>>{}^zZ_J\tag a$$
where 
$$c(h,g)=g,\qua c'(h,g)=({}^hP,{}^hP',hgh\i U_{{}^hP}).$$
Now $c$ is a smooth morphism with fibres isomorphic to $G$ and $c'$ is a smooth morphism with fibres 
isomorphic to $P\cap P'$. By 2.1 (for $\hL,\hL^1$ instead of $\hG,G^1$) the notion of unipotent character
sheaf on $\hL^1$ is well defined. If $A_1\in CS(\hL^1)$ then $c^\bst A_1$ is a simple 
$(P\cap P')$-equivariant perverse sheaf on $G\T\hL^1$ (for the free $P\cap P'$ action 
$a:(h,g)\m(ha\i,\text{\rm pr}(a)g\text{\rm pr}(a)\i)$ where $\text{\rm pr}:P\cap P'@>>>L$ is the canonical 
projection) hence it is of the form $c'{}^\bst A$ for a well defined
(necessarily $G$-equivariant) simple perverse sheaf $A$ on ${}^zZ_J$. By definition, $CS({}^zZ_J)$ consists
of the simple perverse sheaves $A$ obtained
 as above from some $A_1\in CS(\hL^1)$. Note that $A_1\m A$ defines
a bijection between the set of isomorphism classes of objects in $CS(\hL^1)$ and the set of isomorphism 
classes of objects in $CS({}^zZ_J)$. This definition of $CS({}^zZ_J)$ does not depend on the choice of 
$(P,P')$ in ${}^zY_J$ and that of $L$ in $\cs_{P,P'}$ since the set of triples $(P,P',L)$ as above is a 
homogeneous $G$-space.

Next we assume that $z\n N_{J,\d}$. Then, setting $J_1=J\cap\d\i(\Ad(z)J)$, we have $\sha J_1<\sha J$ so
that $CS({}^wZ_{J_1})$ is defined by the induction hypothesis. By definition, $CS({}^wZ_J)$ consists of the
simple perverse sheaves of the form $\vt_{J,w}^\bst A$ where $A\in CS({}^wZ_{J_1})$ and $\vt_{J,w}$ is as
in 1.6(b). Note that $A\m\vt_{J,w}^\bst A$ defines a bijection from the set of isomorphism classes of 
objects in $CS({}^wZ_{J_1})$ and the set of isomorphism classes of objects in $CS({}^wZ_J)$. 

This completes the inductive definition of $CS({}^wZ_J)$. Note that if $A\in CS({}^wZ_J)$ then $A$ is
$G$-equivariant. Note that $\bbq\la|{}^wZ_J|\ra\in CS({}^wZ_J)$.

\subhead 2.3\endsubhead
Let $J\sub I$ and let $w\in{}^{{}^\d J}W$. Let $A\in CS({}^wZ_J)$. Let $A^\sha$ be the unique simple 
perverse 
sheaf on $Z_J$ such that $A^\sha|_{{}^wZ_J}=A$ and $\supp(A^\sha)$ is the closure in $Z_J$ of $\supp(A)$. 
(Here $\supp$ denotes support.) Let $CS'(Z_J)$ be the collection of simple perverse sheaves on $Z_J$ that 
are isomorphic to $A^\sha$ for some $w,A$ as above. We show that, if $\tA\in CS'(Z_J)$ and $\tA\cong A^\sha$
with $w,A$ as above then 

(a) {\it $w$ is uniquely determined,}

(b) {\it $A$ is uniquely determined up to isomorphism.}
\nl
Assume that we have also $\tA\cong A'{}^\sha$ where $A'\in CS({}^{w'}Z_J)$. Then $\supp(A)$ and 
$\supp(A')$ are dense in $\supp(\tA)$. Hence they have nonempty intersection. Since $\supp(A)\sub{}^wZ_J$,
$\supp(A')\sub{}^{w'}Z_J$, it follows that ${}^wZ_J\cap{}^{w'}Z_J\ne\emp$ so that $w=w'$. We have 
$\tA|_{{}^wZ_J}\cong A$, $\tA|_{{}^wZ_J}\cong A'$ hence $A\cong A'$.

We now state the following result.

\proclaim{Proposition 2.4} Let $J\sub I$. 

(a) We have $CS(Z_J)=CS'(Z_J)$.

(b) Let $w\in{}^{{}^\d J}W$ and let $A\in CS(Z_J)$. Let $A'$ be a simple perverse sheaf on ${}^wZ_J$ such 
that $A'\dsv A|_{{}^wZ_J}$. Then $A'\in CS({}^wZ_J)$.
\endproclaim
The proof of (a) appears in \cite{\PCS, 4.13, 4.17}. The proof of (b) appears in \cite{\PCS, 4.12}. We will 
reprove them here (the proof of 2.4(b) is given in 2.12; the proof of 2.4(a) is given in 2.14). To do so
we are using a number of lemmas some of which are more precise than those in \cite{\PCS}.

\proclaim{Lemma 2.5}Let $J\sub I$, $z\in N_{J,\d}$ (so that $z\in{}^{{}^\d J}W$). Let $y\in W$ and let $A$ 
be a simple perverse sheaf on ${}^zZ_J$ such that $A\dsv\KK^y_J|_{{}^zZ_J}$. Then $A\in CS({}^zZ_J)$.
\endproclaim
If $y\n z\i W_{{}^\d J}$ then the image of $\p_{J,y}:\cb^2_{J,y}@>>>Z_J$ is disjoint from ${}^zZ_J$ hence 
$\KK^y_J|_{{}^zZ_J}=0$. This contradicts the choice of $A$. Thus we must have $y=z\i y'$ for some 
$y'\in W_{{}^\d J}$. We fix $(P,P')\in{}^zY_J$, $L\in\cs_{P,P'}$ (see 2.2). Let $\cb_L$ be the variety of 
Borel subgroups of $L$. Let $\cb_{L,y'}=\{(\b,\b')\in\cb_L^2;(\b U_{P'},\b'U_{P'})\in\co_{y'}$. We have a 
cartesian diagram
$$\CD \Xi@<\tc<< G\T\Xi@>\tc'>>\cb^2_{J,z\i y'}\\
  @VfVV      @Vf'VV        @Vf''VV        \\
\hL^1@<c<<G\T\hL^1@>c'>>{}^zZ_J
\endCD$$
where 
$$\Xi=\{(\b,\b',g)\in\cb_L\T\cb_L\T\hL^1;g\bg\i=\b',(\b,\b')\in\co_{L,y'}\},$$
$$\tc(h,(\b,\b',g))=(\b,\b',g),$$ 
$$\tc'(h,(\b,\b',g))=(hP^{\b U_{P'}}h\i,h\b'U_{P'}h\i,hgh\i U_{hPh\i}),$$
$$f(\b,\b',g)=g,\qua f'(h,(\b,\b',g))=(h,g),\qua f''(B,B',gU_{B_J})=(B_J,B'_{{}^\d J},gU_{B_J})$$
and $c,c'$ are as in 2.2(a). It follows that $c^*f_!\bbq=c'{}^*f''_!\bbq$. (Both are equal to $f'_!\bbq$.)
We have $f_!\bbq=\KK^{y';\hL^1}$ (which is defined like $\KK^y_J$ by replacing $\hG,G^1,J,y$ by
$\hL,\hL^1,{}^\d J,y'$) and $f''_!\bbq=\KK^{z\i y'}_J|_{{}^zZ_J}$. Thus we have
$$c'{}^*(\KK^{z\i y'}_J|_{{}^zZ_J})=c^*(\KK^{y';\hL^1}).\tag a$$
Since $c$ is smooth with fibres isomorphic to $G$ and $c'$ is smooth with fibres isomorphic to $P\cap P'$ it
follows that
$$c'{}^\bst\op_j{}^pH^j(\KK^{z\i y'}_J|_{{}^zZ_J})=c^\bst\op_j{}^pH^j(\KK^{y';\hL^1})\tag b$$
so that $c'{}^\bst A\dsv c^\bst\KK^{y';\hL^1}$ hence there exists a simple perverse sheaf $C$ on 
$\hL^1$ such that $c'{}^\bst A=c^\bst C$ and $C\dsv\KK^{y';\hL^1}$. Thus $C\in CS(\hL^1)$ and 
from the definitions we see that $A\in CS({}^zZ_J)$. The lemma is proved.

\subhead 2.6\endsubhead
For any $y_1,y_2,y_3$ in $W$ and any $i\in\ZZ$ we set
$$R^i_{y_1,y_2,y_3}=H^i_c(\{B\in\cb;(B_1,B)\in\co_{y_1},(B,B'_1)\in\co_{y_2}\},\bbq)$$
where $(B_1,B'_1)\in\co_{y_3}$ is fixed. This is a $\bbq$-vector space independent of the choice of 
$(B_1,B'_1)$ (since $G$ acts on $\co_{y_3}$ transitively with connected isotropy groups). 

For $J\sub I$, $u\in W$ we define $\tp_{J,u}:{}^{u\i}\tZ_\emp@>>>\tZ_J$ by $(B,B',g)\m(B_J,B'_{{}^\d J},g)$ 
and we set $\tKK^u_J=(\tp_{J,u})_!\bbq\in\cd(\tZ_J)$; we have $\tKK^u_J=e_J^*\KK^u_J$. 

\proclaim{Lemma 2.7} Let $J\sub I$. Let $y^*\in W^{{}^\d J}$, $y_*\in W_{{}^\d J}$, $y=y^*y_*$. 
Let $z\in{}^{{}^\d J}W^J$ and let $\Om=W_{{}^\d J}zW_J$. We set $J_1=J\cap\d\i(\Ad(z)J)$. Let 
$\ti\a:\tZ_{J,\Om}@>\si>>\tZ^\da_{J_1,\Om}$ be as in 1.5(b). 

(a) In $\cd(\tZ_{J,\Om})$ we have
$$\tKK^y_J|_{\tZ_{J,\Om}}\Bpq
\{R^i_{\d\i(y_*),y^*;y'y^*}\ot\ti\a^*(\tKK^{y'y^*}_{J_1}|_{\tZ^\da_{J_1,\Om}})[-i];y'\in W_J,i\in\ZZ\}.$$

(b) Let $w\in\Om\cap{}^{{}^\d J}W$. Let $\vt_{J,w}:{}^wZ_J@>\si>>{}^wZ_{J_1}$ be as in 1.6(b). In 
$\cd({}^wZ_J)$ we have
$$\KK^y_J|_{{}^wZ_J}\Bpq\{R^i_{\d\i(y_*),y^*;y'y^*}
\ot\vt_{J,w}^*(\KK^{y'y^*}_{J_1}|_{{}^wZ_{J_1}})[-i];y'\in W_J,i\in\ZZ\}.$$

(c) In the setup of (b) let $\t_w:W_{{}^\d(J^w_\iy)}@>>>W_{{}^\d(J^w_\iy)}$ be as in 1.3. Let 
$x\in W_{{}^\d(J^w_\iy)}$. We have
$$\KK^{w\i x}_J|_{{}^wZ_J}=\vt_{J,w}^*(\KK^{w\i\t_w(x)}_{J_1}|_{{}^wZ_{J_1}}).$$
\endproclaim
We prove (a). Assume first that $y\i\n\Om$. In this case, the image of $\tp_{J,y}:{}^{y\i}\tZ_\emp@>>>\tZ_J$
is disjoint from $\tZ_{J,\Om}$ hence $\tKK^y_J|_{\tZ_{J,\Om}}=0$. Moreover, for any $y'\in W_J$, the image 
of $\tp_{J_1,y'y^*}:{}^{(y'y^*)\i}\tZ_\emp@>>>\tZ_{J_1}$ is disjoint from $\tZ^\da_{J_1,\Om}$. (We have 
$y^{*-1}\n\Om$. If $y^{*-1}y'{}\i\in W_{{}^\d J_1}zW_J$ then $y^{*-1}\in W_{{}^\d J_1}zW_J\sub\Om$, a 
contradiction.) Thus $\tKK^{y'y^*}_{J_1}|_{\tZ^\da_{J_1,\Om}}=0$ and (a) holds.

We now assume that $y\i\in\Om$. We set
$$\align&E^y_J=\\&\{(B_1,B,B'_1,g)\in\cb\T\cb\T\cb\T G^1;
(B_1,B)\in\co_{\d\i(y_*)},(B,B'_1)\in\co_{y^*},B'_1={}^gB_1\},\endalign$$
$${}^{y^{*-1}W_J}\tZ_\emp=\cup_{y'\in W_J}{}^{(y'y^*)\i}\tZ_\emp\sub\tZ_\emp.$$
Note that $\th:E^y_J@>>>{}^{y\i}\tZ_\emp$, $(B_1,B,B'_1,g)\m(B,{}^gB,g)$ is a well defined isomorphism.
Define $k:E^y_J@>>>{}^{y^{*-1}W_J}\tZ_\emp$ by $(B_1,B,B'_1,g)\m(B_1,B'_1,g)$.

Now $\tp_{J,y}:{}^{y\i}Z_\emp@>>>\tZ_J$ factors as ${}^{y\i}Z_\emp@>\ph>>\tZ_{J,\Om}@>j>>\tZ_J$ where $j$ is
the inclusion and ${}^{y^{*-1}W_J}\tZ_\emp@>>>\tZ_{J_1}$ (restriction of $\tp_{J_1}$) factors as 
${}^{y^{*-1}W_J}\tZ_\emp@>\ps>>\tZ^\da_{J_1,\Om}@>j_1>>\tZ_{J_1}$ where $j_1$ is the inclusion. We have a 
commutative diagram
$$\CD 
E^y_J@>\th>> {}^{y\i}Z_\emp@>\ph>>\tZ_{J,\Om}@>j>>\tZ_J\\
 @Vk VV                @.             @V\ti\a VV           @.\\
{}^{y^{*-1}W_J}\tZ_\emp@>\ps>>{}@>>>\tZ^\da_{J_1,\Om}@>j_1>>\tZ_{J_1}
\endCD$$
From the definitions we have
$$\k_!\bbq\Bpq\{R^i_{\d\i(y_*),y^*;y'y^*}\ot\tj_{y'!}\bbq[-i]; y'\in W_J,i\in\ZZ\},$$
where 
$\tj_{y'}:{}^{(y'y^*)\i}\tZ_\emp@>>>{}^{y^{*-1}W_J}\tZ_\emp$ is the obvious inclusion. It follows that
$$\ps_!\k_!\bbq\Bpq\{R^i_{\d\i(y_*),y^*;y'y^*}\ot\ps_!\tj_{y'!}\bbq[-i]; y'\in W_J,i\in\ZZ\},$$
that is,
$$\ti\a_!\ph_!\th_!\bbq\Bpq\{R^i_{\d\i(y_*),y^*;y'y^*}\ot\ps_!\tj_{y'!}\bbq[-i]; y'\in W_J,i\in\ZZ\}.$$
Since $\th$ is an isomorphism we have $\th_!\bbq=\bbq$; since $\ti\a$ is an isomorphism we have 
$\ti\a^*\ti\a_!=1$ hence
$$\ph_!\bbq\Bpq\{R^i_{\d\i(y_*),y^*;y'y^*}\ot\ti\a^*\ps_!\tj_{y'!}\bbq[-i]; y'\in W_J,i\in\ZZ\}.$$
Since $j^*j_!=1$, $j_1^*j_{1!}=1$ and $j_!\ph_!\bbq=(\tp_{J,y})_!\bbq=\ti\KK^y_J$,
$$j_{1!}\ps_!\tj_{y'!}\bbq=(\tp_{J_1,y'y^*})_!\bbq=\tKK^{y'y^*}_{J_1},$$
we have
$$j^*\tKK^y_J\Bpq\{R^i_{\d\i(y_*),y^*;y'y^*}\ot\ti\a^*j_1^*(\tKK^{y'y^*}_{J_1})[-i];y'\in W_J,i\in\ZZ\}$$
and (a) is proved.

We prove (b). We have a commutative diagram
$$\CD
{}^w\tZ_J@>h>>\tZ_{J,\Om}\\
@V\vt_{J,w}VV     @V\a VV\\
{}^w\tZ_{J_1}@>h_1>>\tZ^\da_{J_1,\Om}
\endCD$$
where $h,h_1$ are the inclusions. Applying $h^*$ to the relation $\Bpq$ in (a) we obtain
$$\tKK^y_J|_{{}^w\tZ_J}\Bpq\{R^i_{\d\i(y_*),y^*;y'y^*}\ot h^*\ti\a^*j_1^*(\tKK^{y'y^*}_{J_1})[-i];
y'\in W_J,i\in\ZZ\}.$$
Let $\ti\vt_{J,w}$ be as in 1.6(b). Using that $h^*\ti\a^*j_1^*=\ti\vt_{J,w}^*h_1^*j_1^*$ and that 
$j_1h_1:{}^w\tZ_{J_1}@>>>\tZ_{J_1}$ is the inclusion we obtain
$$\tKK^y_J|_{{}^w\tZ_J}\Bpq
\{R^i_{\d\i(y_*),y^*;y'y^*}\ot\ti\vt_{J,w}^*(\tKK^{y'y^*}_{J_1}|_{{}^w\tZ_{J_1}})[-i];y'\in W_J,i\in\ZZ\}$$
or equivalently
$$(e_J^*\KK^y_J)|_{{}^w\tZ_J}\Bpq
\{R^i_{\d\i(y_*),y^*;y'y^*}\ot\ti\vt_{J,w}^*((e_{J_1}^*\KK^{y'y^*}_{J_1})|_{{}^w\tZ_{J_1}})[-i];y'\in W_J,
i\in\ZZ\}.$$
Let $e_{J,w}:{}^w\tZ_J@>>>{}^wZ_J$ (resp. $e_{J_1,w}:{}^w\tZ_{J_1}@>>>{}^wZ_{J_1}$) be the restiction of 
$e_J$ (resp. $e_{J_1}$). We have $(e_J^*\KK^y_J)|_{{}^w\tZ_J}=e_{J,w}^*(\KK^y_J|_{{}^wZ_J})$ and
$$\ti\vt_{J,w}^*((e_{J_1}^*\KK^{y'y^*}_{J_1})|_{{}^w\tZ_{J_1}})=
\ti\vt_{J,w}^*((e_{J_1,w}^*(\KK^{y'y^*}_{J_1}|_{{}^wZ_{J_1}}))=
e_{J,w}^*\vt_{J,w}^*(\KK^{y'y^*}_{J_1}|_{{}^wZ_{J_1}})$$
hence
$$e_{J,w}^*(\KK^y_J|_{{}^wZ_J})\Bpq\{R^i_{\d\i(y_*),y^*;y'y^*}\ot 
e_{J,w}^*\vt_{J,w}^*(\KK^{y'y^*}_{J_1}|_{{}^wZ_{J_1}})[-i];y'\in W_J,i\in\ZZ\}.$$
Applying $(e_{J,w})_!$ and using that $(e_{J,w})_!e_{J,w}^*=\la-2|\cp_J|\ra$ we see that (b) holds.

We prove (c). We shall use notation in the proof of (a) with $y^*=w\i$, $y_*=x$, $y=w\i x$. We have 
$\d\i(x)w\i=w\i\t_w(x)$. Since $w\i\in W^{{}^\d J}$ and $\t_w(x)\in W_{{}^\d J}$ we have 
$$\align&l(\d\i(x)w\i)=l(w\i\t_w(x))=l(w\i)+l(\t_w(x))=l(w\i)+l(x)\\&=l(\d\i(x))+l(w\i).\endalign$$
Hence if $(B_1,B,B'_1,g)\in E^y_J$ then $B$ is uniquely determined by $B_1,B'_1$ and $k$ is an isomorphism 
of $E^y_J$ onto the subspace 
$${}^{w\d\i(x)\i}\tZ_\emp={}^{\t_w(x)\i w}\tZ_\emp$$ 
of ${}^{y^{*-1}W_J}\tZ_\emp$. Hence we have 
$j_{1!}\ti\a_!\ph_!\th_!\bbq=j_{1!}\ps_!k_!\bbq=\tKK^{w\i\t_w(x)}_{J_1}$. Since $\th$ is an isomorphism we 
have $\th_!\bbq=\bbq$. Hence 
$j_{1!}\ti\a_!\ph_!\bbq=\tKK^{w\i\t_w(x)}_{J_1}$.
Since $j_1^*j_{1!}=1$, $\ti\a^*\ti\a_!=1$ we deduce $\ph_!\bbq=\ti\a^*j_1^*\tKK^{w\i\t_w(x)}_{J_1}$. We have
$j_!\ph_!\bbq=\tKK^{w\i x}_J$. Since $j^*j_!=1$, we deduce $\ph_!\bbq=j^*\tKK^{w\i x}_J$. Thus
$$j^*\tKK^{w\i x}_J=\ti\a^*j_1^*\tKK^{w\i\t_w(x)}_{J_1}.$$
Applying $h^*$ (notation in the proof of (b)) we obtain
$$\tKK^{w\i x}_J|_{{}^w\tZ_J}=\ti\vt_{J,w}^*(\tKK^{w\i\t_w(x)}_{J_1}|_{{}^w\tZ_{J_1}}).$$
From this we deduce as in the proof of (b) that (c) holds. The lemma is proved.

\proclaim{Lemma 2.8}Let $J\sub I$, $w\in{}^{{}^\d J}W$. Let $y\in W$. Let $A$ be a simple perverse sheaf on
${}^wZ_J$ such that $A\dsv(\KK_y|_{{}^wZ_J})$. Then $A\in CS({}^wZ_J)$.
\endproclaim
We argue by induction on $\sha J$. Let $z=\min(W_{{}^\d J}wW_J)$.

Assume first that $z\in N_{J,\d}$ so that $z=w$. If $y\n z\i W_{{}^\d J}$ then the image of 
$\p_{J,y}:\cb^2_{J,y}@>>>Z_J$ is disjoint from ${}^wZ_J$ hence $\KK^y|_{{}^wZ_J}=0$. This contradicts the 
choice of $A$. Thus we must have $y=z\i y'$ for some $y'\in W_{{}^\d J}$. Then with the notation in 2.5(b) 
we have
$$c'{}^\bst\op_j{}^pH^j(\KK^{z\i y'}|_{{}^wZ_J})=c^\bst\op_j{}^pH^j(\KK^{y';\hL^1}).$$
Hence $c'{}^\bst A\dsv c^\bst\op_j{}^pH^j(\KK^{y';\hL^1})$ and there exists a simple perverse sheaf
$A_1$ on $\hL^1$ such that $c'{}^\bst A=c^\bst A_1$ and $A_1\dsv\KK^{y';\hL^1}$. Thus 
$A_1\in CS(\hL^1)$ and from the definitions we see that $A\in CS({}^wZ_J)$.

Next we assume that $z\n N_{J,\d}$. Then, setting $J_1=J\cap\d\i(\Ad(z)J)$, we have $\sha J_1<\sha J$. We 
can write uniquely $y=y^*y_*$ as in 2.7. From 2.7(b) we see that there exists $y'\in W_J$ such that 
$A\dsv\vt_{J,w}^*(\KK^{y'y^*}_{J_1}|_{{}^wZ_{J_1}})$. Hence there exists a simple perverse sheaf $A_1$ on 
${}^wZ_{J_1}$ such that $A=\vt_{J,w}^\bst A_1$ and $A_1\dsv\KK^{y'y^*}_{J_1}|_{{}^wZ_{J_1}}$. By the 
induction hypothesis, we have $A_1\in CS({}^wZ_{J_1})$. Hence $A\in CS({}^wZ_J)$. The lemma is proved.

\proclaim{Lemma 2.9}Let $J\sub I$ and let $w\in{}^{{}^\d J}W$. Let $A\in CS({}^wZ_J)$. 

(a) There exists $x\in W_{{}^\d(J^w_\iy)}$ such that $A\dsv\KK^{w\i x}_J|_{{}^wZ_J}$.

(b) There exists $x\in W_{{}^\d(J^w_\iy)}$ such that $A^\sha\dsv\KK^{w\i x}_J$. In particular, 
$A^\sha\in CS(Z_J)$. Thus we have $CS'(Z_J)\sub CS(Z_J)$.
\endproclaim
We prove (a) by induction on $\sha J$. Let $z=\min(W_{{}^\d J}wW_J)$. Let 
$\t_w:W_{{}^\d(J^w_\iy)}@>>>W_{{}^\d(J^w_\iy)}$ be as in 1.3.

Assume first that $z\in N_{J,\d}$ so that $z=w$ and $J=J^w_\iy$. With notation in 2.2(a) we have 
$c'{}^\bst A=c^\bst A_1$ where $A_1\in CS(\hL^1)$ so that $A_1\dsv\KK^{x;\hL^1}$ for some 
$x\in W_{{}^\d J}$. Since 
$$c'{}^\bst\op_j{}^pH^j(\KK^{w\i x}|_{{}^wZ_J})=c^\bst(\op_j{}^pH^j(\KK^{x;\hL^1}))$$
(see 2.5(b)) we see that $c'{}^\bst A\dsv c'{}^\bst\KK^{w\i x}|_{{}^wZ_J}$ hence 
$A\dsv\KK^{w\i x}|_{{}^wZ_J}$. Thus (a) holds in this case.

Next we assume that $z\n N_{J,\d}$. Then, setting $J_1=J\cap\d\i(\Ad(z)J)$ we have $\sha J_1<\sha J$. We 
have $A=\vt_{J,w}^\bst A_1$ where $A_1\in CS({}^wZ_{J_1})$. By the induction hypothesis, there exists 
$x\in W_{{}^\d(J^w_\iy)}$ such that $A_1\dsv\KK^{w\i x}_{J_1}|_{{}^wZ_{J_1}}$. Hence 
$$A=\vt_{J,w}^*A_1\dsv\vt_{J,w}^*(\KK^{w\i x}_{J_1}|_{{}^wZ_{J_1}}).$$
Using 2.7(c) with $x$ replaced by $\t_w\i(x)$, we deduce that $A\dsv\KK^{w\i\t_w\i(x)}_J|_{{}^wZ_J}$. This 
proves (a).

We prove (b). Let $x$ be as in (a). Applying \cite{\CDGVIII, 36.3(c)} with $Y,Y',C$ replaced by 
${}^wZ_J,Z_J,\KK^{w\i x}_J|_{{}^wZ_J}$ we deduce that (b) holds. (We use the fact that, if $i:{}^wZ_J@>>>Z_J$
is the inclusion, then $i_!(\KK^{w\i x}_J|_{{}^wZ_J})=\KK^{w\i x}_J$, which follows from 1.8(b).)

\subhead 2.10\endsubhead
Let $J\sub I$. For $y\in W$ let $\bco_y$ be the closure of $\co_y$ in $\cb^2$. The closure of $\cb^2_{J,y}$ 
in $\cb^2_J$ is
$$\bcb^2_{J,y}=\{(B,B',gU_{B_J})\in\cb^2_J;(B,B')\in\bco_y\}.$$ 
 
Let $\ck_{J,y}$ be the intersection cohomology complex of $\bcb^2_{J,y}$ with coefficients in $\bbq$, 
extended by zero on $\cb^2_J-\bcb^2_{J,y}$. Let $\bKK^y_J=\p_{J!}\ck_{y,J}$ with $\p_J:\cb^2_J@>>>Z_J$ as in
1.8. For $y'\in W$ let $\io_{y'}:\cb^2_{J,y'}@>>>\cb^2_J$ be the inclusion. By \cite{\KLL} we have
$$\ck_{y,J}\Bpq\{\cv_{y',y}^i\ot\io_{y'!}\bbq\la-i\ra;y'\in W,y'\le y,i\in2\ZZ\}.\tag a$$
where $\cv_{y',y}^i$ are $\bbq$-vector spaces such that
$$\sum_{i\in2\ZZ}\dim\cv_{y',y}^iv^i=P_{y',y}(v^2)\tag b$$
and $P_{y',y}$ is the polynomial in \cite{\KL}. Applying $\p_{J!}$ we obtain
$$\bKK^y_J\Bpq\{\cv_{y',y}^i\ot\KK^{y'}_J\la-i\ra;y'\in W,y'\le y,i\in2\ZZ\}.\tag c$$
Note that
$$|\bcb^2_{J,y}|=l(y)+|G|-|\cp_J|.\tag d$$

\subhead 2.11\endsubhead
Let $J\sub I$ and let $A$ be a simple perverse sheaf on $Z_J$.  We show that conditions (i),(ii) below are 
equivalent:

(i) $A\in CS(Z_J)$;

(ii) {\it We have $A\dsv\bKK^y_J$ for some $y\in W$.}
\nl
Assume that (ii) holds. Then from 2.10(c) we see that $A\dsv\KK^{y'}_J$ for some $y'\in W$, $y'\le y$ hence 
(i) holds. Conversely, assume that (i) holds. We can find $y\in W$ such that $A\dsv\KK^y_J$ for some 
$y\in W$ and such that $A\not\dsv\KK^{y'}_J$ for any $y'\in W$ such that $y'<y$. Using this and 2.10(c) we 
see that (ii) holds. (We use that $\dim \cv_{y,y}^0=1$ and $\cv_{y,y}^i=0$ for $i\ne0$.)

\subhead 2.12\endsubhead
We prove 2.4(b). By 2.11, $A$ is a composition factor of ${}^pH^j(\bKK^y_J)$ for some $y\in W$ and some 
$j\in\ZZ$. By the decomposition theorem, ${}^pH^j(\bKK^y_J)$ is a semisimple perverse sheaf and 
${}^pH^j(\bKK^y_J)[-j]$ is a direct summand of $\bKK^y_J$. It follows that $A[j]$ is a direct summand of 
$\bKK^y_J$. Hence $A'\dsv\bKK^y_J|_{{}^wZ_J}$. Using 2.10(c) we deduce that $A'\dsv\KK^{y'}_J|_{{}^wZ_J}$ 
for some $y'\le y$. Using 2.8 we deduce that $A'\in CS({}^wZ_J)$. This proves 2.4(b).

\proclaim{Lemma 2.13} Let $J\sub I$. Let $A\in CS(Z_J)$. There exists $w\in{}^{{}^\d J}W$ and 
$A'\in CS({}^wZ_J)$ such that $A=A'{}^\sha$. In particular, $CS(Z_J)\sub CS'(Z_J)$.
\endproclaim
The subsets $\{\supp(A)\cap{}^wZ_J;w\in{}^{{}^\d J}W\}$ form a partition of $\supp(A)$ into locally closed
subvarieties. Hence we can find $w\in{}^{{}^\d J}W$ such that $\supp(A)\cap{}^wZ_J$ is open dense in 
$\supp(A)$. We set $A'=A|_{{}^wZ_J}$. Then $A'$ is a simple perverse sheaf on ${}^wZ_J$ and 
$A'\in CS({}^wZ_J)$ (see 2.4(b)). From the definitions we have $A=A'{}^\sha$. The lemma is proved.

\subhead 2.14\endsubhead
Let $J\sub I$. Since $CS'(Z_J)\sub CS(Z_J)$ (see 2.9(b)) and $CS(Z_J)\sub CS'(Z_J)$ (see 2.13) we see that
2.4(a) holds.

\subhead 2.15\endsubhead
For any $J\sub I$ and any $w\in{}^{d(J)}W$ we choose $(P,P')\in\cp_{J^w_\iy}\T\cp_{{}^\d(J^w_\iy)}$ and a 
common Levi subgroup $L^w$ of $P,P'$. Define $\hL^w,\hL^{w1}$ as $\hL,\hL^1$ in 2.2 with $J,L$ replaced by 
$J^w_\iy,L^w$. We shall denote by $\uCS(\hL^{w1})$ a set of representatives for the objects in 
$CS(\hL^{w1})$. We can assume that $\bbq\la|L^w|\ra\in\uCS(\hL^{w1})$.

For $C\in\uCS(\hL^{w1})$ we denote by $\tC$ the object of $CS({}^wZ_{J^w_\iy})$ such that 
$c'{}^\bst\tC=c^\bst C$ (where $c,c'$ are as in 2.5 with $J$ replaced by $J^w_\iy$); we set 
$$C_w=\vt_{J,w}^\bst\vt_{J_1,w}^\bst\vt_{J_2,w}^\bst\do(\tC)\in CS({}^wZ_J)$$ 
where $(J_n,w_n)_{n\ge0}=\k(w)$. Then $C_w^\sha:=(C_w)^\sha$ is defined as in 2.3; note that
$C_w^\sha\in CS(Z_J)$ by 2.4(a). From 2.4 and the definitions we see that  

(a) {\it $\{C_w;C\in\uCS(\hL^{w1})\}$ is a set of representatives for the isomorphism classes of objects in 
$CS({}^wZ_J)$;}

(b) {\it $\{C_w^\sha; w\in{}^{{}^\d J}W,C\in\uCS(\hL^{w1})\}$ is a set of representatives for the 
isomorphism classes of objects in $CS(Z_J)$.}

Let $\fK_J$ be the free $\ca$-module ($\ca$ as in 0.1) with basis
$\{C_w^\sha;w\in{}^{{}^\d J}W,C\in\uCS(\hL^{w1})\}$. For $w\in{}^{{}^\d J}W$ let ${}^w\fK_J$ be the free 
$\ca$-module with basis $\{C_w;C\in\uCS(\hL^{w1})\}$.
Let $\fK(\hL^{w1})$ be the free $\ca$-module with basis $\{C;C\in\uCS(\hL^{w1})\}$.

Let $\x\m\x_w$ be the $\ca$-module isomorphism $\fK(\hL^{w1})@>\si>>{}^w\fK_J$ such that $C\m C_w$ for any
$C\in\uCS(\hL^{w1})$.

\head 3. Mixed structures\endhead
\subhead 3.1\endsubhead
In this section we assume that $\kk$ is an algebraic closure of a finite field $\FF_q$ with $q$ elements and
that we are given an $\FF_q$-rational structure on $\hG$ with Frobenius map $F:\hG@>>>\hG$ such that
$F(G^1)=G^1$ and the restriction of $F$ to $G$ is the Frobenius map of an $\FF_q$-split rational structure 
on $G$. Then for any $J\sub I$, we have $P\in\cp_J\imp F(P)\in\cp_J$ and $P\m F(P)$ is the Frobenius map of
an $\FF_q$-rational structure on $\cp_J$; moreover, $(P,P',gU_P)\m(F(P),F(P'),F(g)U_{F(P)})$ is the 
Frobenius map of an $\FF_q$-rational structure on $Z_J$. For any $w\in W$ we have 
$(B,B')\in\co_w\imp(F(B),F(B'))\in\co_w$ and $(B,B')\m(F(B),F(B'))$ is the Frobenius map of an 
$\FF_q$-rational structure on $\co_w$. For any $J\sub I$ and any $w\in{}^{{}^\d J}W$, ${}^wZ_J$ is a 
subvariety
of $Z_J$ defined over $\FF_q$; we choose $P,P',L^w$ as in 2.15 in such a way that $F(P)=P$, $F(P')=P'$, 
$F(L^w)=L^w$, $F(\hL^{w1})=\hL^{w1}$ (notation of 2.2 with $J$ replaced by $J^w_\iy$). We shall assume (as 
we may) that for any $J\sub I$, any $w\in{}^{{}^\d J}W$, and any $C\in\uCS(\hL^{w1})$ we can find an 
isomorphism $\ph_C:F^*C@>\si>>C$ which makes $C$ into a pure complex of weight $0$; we shall assume that 
such a $\ph_C$ has been chosen.

Let $J\sub I$, $w\in{}^{{}^\d J}W$. For $C\in\uCS(\hL^{w1})$ let $\tC\in CS({}^wZ_{J^w_\iy})$,
$C_w\in CS({}^wZ_J)$, $C_w^\sha\in CS(Z_J)$ be as in 2.15. Then $\tC,C_w,C_w^\sha$ inherit from $C$ mixed 
structures which are pure of weight $0$. 

\subhead 3.2\endsubhead
Let $J\sub I$. Let $\o$ (resp. $\o_J$) be the element of maximal length in $W$ (resp. $W_J$). Let 
$\tB\in\cb$. Let $\tP=\tB_{\o J\o}$, $\tU=U_{\tP}$. Now $\tU$ acts by conjugation on
$$\cu=\{B\in\cb;\pos(\tB,B)\in\o W_J\},$$ 
an open subset of $\cb$. For any $P\in\cp_J$ such that $\pos(\tB,P)=\o\o_J$ we define $\cb_P\T\tU@>>>\cu$ by
$(B_1,u)\m uB_1u\i$. We show:

(a) {\it This map is a bijection (in fact an isomorphism).}
\nl
Assume that $B'\in\cb_P,B''\in\cb_P$ and $u'\in\tU,u''\in\tU$ satisfy $u'B'u'{}\i=u''B''u''{}\i$. Setting 
$u=u''{}\i u'\in\tU$ we have $uB'u\i=B''$ hence $uPu\i=P$ and $u\in P$. Now $\tP,P$ are opposed parabolic 
subgroups so that $P\cap\tU=\{1\}$. Thus $u=1$ so that $u'=u''$ and $B'=B''$. We see that our map is 
injective. Let $B\in\cu$. We have $B\in\cb_{P'}$ where $P'\in\cp_J,\pos(\tB,P')=\o\o_J$. Now $\tU$ acts 
transitively (by conjugation) on $\{P_1\in\cp_J,\pos(\tB,P_1)=\o\o_J\}$ hence there exists $u\in\tU$ such 
that $uPu\i=P'$. Setting $B_1=u\i Bu$ we have $B_1\in\cb_P$ and our map takes $(B_1,u)$ to $B$; thus it is
surjective hence bijective. 

Let $L=\tP\cap P$; this is a common Levi subgroup of $\tP$ and $P$. Let $S$ be the identity component of the
centre of $L$. We can find a one parameter subgroup $\l:\kk^*@>>>S$ such that $\lim_{t\m0}\l(t)u\l(t)\i=1$ 
for any $u\in\tU$. We define an action $t:B\m\l(t)B\l(t)\i$ of $\kk^*$ on $\cu$. (This is well defined since
$\l(t)\in\tB$.) Under the isomorphism (a) this action becomes the action of $\kk^*$ on $\cb_P\T\tU$ given by
$t:(B_1,u)\m(B_1,\l(t)u\l(t)\i)$. (To see this we must check that 
$\l(t)u\l(t)\i B_1\l(t)u\i\l(t)\i=\l(t)uB_1u\i\l(t)\i$ for $B_1\in\cb_P$; this holds since $\l(t)\in B_1$.)
Now $\lim_{t\m0}(B_1,\l(t)u\l(t)\i)=(B_1,1)$ for any $(B_1,u)\in\cb_P\T\tU$. Hence $t:B\m\l(t)B\l(t)\i$ is a
flow on $\cu$ which contracts $\cu$ to its fixed point set $\cb_P$. We have 
$\cb\T\cp_J=\sqc_{y\in W^J}\co'_y$ where $\co'_y$ is the image of $\cb^2_{yW_J}:=\cup_{a\in W_J}\co_{ya}$ 
under the obvious map $\cb\T\cb@>>>\cb\T\cp_J$. This is exactly the decomposition of $\cb\T\cp_J$ into 
$G$-orbits where $G$ acts by simultaneous conjugation. Hence $\cb^2_{yW_J}$ is a locally closed subvariety 
of $\cb^2$ for any $y\in W^J$.

\subhead 3.3\endsubhead
Let $J\sub I$.
We fix $y\in W^J$. Let $\tB\sub\tP$ be as in 3.2. Let $B_1\in\cb$ be such that $(\tB,B_1)\in\co_\o$ and 
let $P=(B_1)_J$; then $\pos(\tB,P)=\o\o_J$. We define $S,\l$ in terms of $\tP,P$ as in 3.2. We choose 
$B^*\in\cb$ as follows: we note that $T_1:=\tB\cap B_1$ is a maximal torus of $G$ containing $S$ (since 
$S\sub\tB$, $S\sub B_1$) and we choose $B^*$ so that $T_1\sub B^*$ and $(\tB,B^*)\in\co_{\o y\i}$. We have 
$S\sub B^*$. Since $(B^*,\tB)\in\co_{y\o}$, $(\tB,B_1)\in\co_\o$ and $B^*,\tB,B_1$ contain a common 
maximal torus, we have $(B^*,B_1)\in\co_{y\o\o}=\co_y$. Hence for any $B\in\cb_P$ we have 
$(B^*,B)\in\co_{ya}$ with $a\in W_J$. In other words, we have $\cb_P\sub\cup_{a\in W_J}\cb_{ya}$, where for 
any $z\in W$ we set $\cb_z=\{B\in\cb;(B^*,B)\in\co_z\}$. As in 3.2, we set 
$\cu=\{B\in\cb;\pos(\tB,B)\in\o W_J\}$. Note that $\cup_{a\in W_J}\cb_{ya}\sub\cu$. (If $B\in\cb$ satisfies
$\pos(B^*,B)\in yW_J$ then for some $B'\in\cb$ we have $\pos(B^*,B')=y,\pos(B',B)\in W_J$ hence 
$\pos(\tB,B')=\o$ and $\pos(\tB,B)\in\o W_J$.)

For $z\in W$ we set $\bcb_z=\{B\in\cb;(B^*,B)\in\bco_z\}$. Let $w\in W$. We show:

(a) {\it $\bKK^w|_{\cb^2_{yW_J}}$ is pure of weight $0$.}
\nl
Let $\bck^w$ be the intersection cohomology complex of $\bcb_w$ with coefficients in $\bbq$; it is naturally
a pure complex of weight $0$. Using the transitivity of the simultaneous conjugation action of $G$ on 
$\co'_y$ and the fact the fibre of $\cb^2_{yW_J}@>>>\co'_y$ at $(B^*,P)\in\co'_y$ may be identified with
$\cb_P$, we see that it is enough to show that:

(b) {\it $\bck^w|_{\bcb_w\cap\cb_P}$ is pure of weight $0$.}
\nl
We can assume that $\bcb_w\cap\cb_P\ne\emp$. Since $\bcb_w\cap\cu$ is open in $\bcb_w$, we have 
$\bck^w|_{\bcb_w\cap\cu}=K$ where $K$ is the intersection cohomology complex of $\bcb_w\cap\cu$ and it is 
enough to show that $K|_{\bcb_w\cap\cb_P}$ is pure of weight $0$ (recall that $\cb_P\sub\cu$). For any 
$z\in W$, $\cb_z\cap\cu$ is stable under the $\kk^*$-action $t:B\m\l(t)B\l(t)\i$ on $\cu$ (we use that 
$\l(t)\in B^*$ for any $t$). Hence $\bcb_w\cap\cu$ is stable under this $\kk^*$-action on $\cu$. Since the 
$\kk^*$-action on $\cu$ is a contraction to its fixed point set $\cb_P$ and $\bcb_w\cap\cu$ is closed in 
$\cu$ and $\kk^*$-stable, we deduce that the $\kk^*$-action on $\bcb_w\cap\cu$ is a contraction to 
$\bcb_w\cap\cb_P$ so that (b) follows from the ``hyperbolic localization theorem'' \cite{\BRA}. This proves 
(a).

\proclaim{Proposition 3.4}Let $J\sub I$ and let $z\in N_{J,\d}$. 

(a) For $y\in W$, $\bKK^y_J|_{{}^zZ_J}$ (with its natural mixed structure) is pure of weight $0$.

(b) If $w\in{}^{{}^\d J}W$, $C\in CS(\hL^{w1})$, then 
$C_w^\sha|_{{}^zZ_J}$ (with its natural mixed structure) is pure of weight $0$.
\endproclaim
We prove (a). We have ${}^zZ_J=\{(P,P',gU_P)\in Z_J;\pos(P',P)=z\}$. We have a diagram
$$\cb^2_{z\i W_{{}^\d J}}@<c<<\Xi@>d>>{}^zZ_J$$
where $\Xi=\{(B,B',gU_{B_J})\in\cb^2_J; (B,B')\in\cup_{a\in W_{{}^\d J}}\co_{z\i a}\}$
and $c(B,B',gU_{B_J})=(B,B')$, $d(B,B',gU_{B_J})=(B_J,B'_{{}^\d J},gU_{B_J})$. 
Now $\Xi$ is the inverse image of ${}^zZ_J$ under $\p_J:\cb^2_J@>>>Z_J$ and $d$ is the restriction of $\p_J$.
Moreover $\p_J$ is proper hence $d$ is proper. Note also that $c$ is smooth. From the definitions we see that
$\bKK^y_J|_{{}^zZ_J}=d_!c^*(\bKK^y|_{\cb^2_{z\i W_{{}^\d J}}})$. It remains to note that 
$\bKK^y|_{\cb^2_{z\i W_{{}^\d J}}}$ is pure of weight $0$ (see 3.3(a) with $J$ replaced by ${}^\d J$), that 
$c^*$ maps pure complexes of weight zero to pure complexes of weight zero (since $c$ is smooth) and $d_!$ 
maps pure complexes of weight zero to pure complexes of weight zero (by Deligne's theorem applied to the 
proper map $d$).

We prove (b). We can find $y\in W$ and $j\in\ZZ$ such that $C_w^\sha$ appears in ${}^pH^j(\bKK^y_J)$ (with
mixed structures being not necessarily compatible). Since $\bKK^y_J$ and $\bKK^y_J|_{{}^zZ_J}$ are pure of 
weight $0$ it follows that ${}^pH^j(\bKK^y_J)$ and ${}^pH^j(\bKK^y_J)|_{{}^zZ_J}$ are pure of weight $j$. We
can find a nonzero mixed vector space $V$ of pure weight $j$ such that $V\ot C_w^\sha$ is a direct 
summand of ${}^pH^j(\bKK^y_J)$ (respecting the mixed structures). Then $V\ot C_w^\sha|_{{}^zZ_J}$ is a 
direct summand of ${}^pH^j(\bKK^y_J)|_{{}^zZ_J}$ (respecting the mixed structures). Hence 
$C_w^\sha|_{{}^zZ_J}$ is pure of weight $0$.

{\it Remark.} More generally for $J\sub I$, $y\in W$, $w\in{}^{{}^\d J}W$, we expect that 
$\bKK^y_J|_{{}^wZ_J}$ (with its natural mixed structure) is pure of weight $0$.

\subhead 3.5\endsubhead
Let $J\sub I$. Let $y\in W$. Since $\bKK^y_J\la|\bcb^2_{J,y}|\ra$ is pure of weight $0$, we have for any 
$j\in\ZZ$:
$${}^pH^{-j}(\bKK^y_J\la|\bcb^2_{J,y}|\ra)=
\op_{w\in{}^{{}^\d J}W,C\in\uCS(\hL^{w1})}V_{y,w,C,j}(j/2)\ot C_w^\sha \tag a$$
in $\cd_m(Z_J)$, where $V_{y,w,C,j}$ are $\bbq$-vector spaces of pure weight $0$. Moreover, by the relative 
hard Lefschetz theorem \cite{\BBD, 5.4.10}, for any $y,w,C,j$ we have:
$$V_{y,w,C,j}=V_{y,w,C,-j}.\tag b$$
Hence if $t\in N_{J,\d}$ we have
$$\bKK^y_J|_{{}^tZ_J}\la|\bcb^2_{J,y}|\ra\Bpq  
\{V_{y,w,C,j}\ot C_w^\sha|_{{}^tZ_J}\la j\ra;w\in{}^{{}^\d J}W,C\in\uCS(\hL^{w1}),j\in\ZZ\}\tag c$$
Since for $w\in{}^{{}^\d J}W$, $C_w^\sha|_{{}^tZ_J}$ is pure of weight $0$ (see 3.4) we have 
$$C_w^\sha|_{{}^tZ_J}\Bpq\{{}'V{}_{t,w,j'}^{C',C}\ot C'_t\la j'\ra; C'\in\uCS(\hL^{t1}),j'\in\NN\}\tag d$$
in $\cd_m({}^tZ_J)$, where ${}'V_{t,w,j'}^{C',C}$ are $\bbq$-vector spaces of pure weight $0$.
Note that if $w=t\in N_{J,\d}$ we have $C_w^\sha|_{{}^tZ_J}=C_t$ 
hence ${}'V_{t,t,j'}^{C',C}$ is $\bbq$ if $C'=C,j'=0$ and ${}'V_{t,t,j'}^{C',C}$ is $0$ otherwise.
From (c),(d) we deduce
$$\align&\bKK^y_J|_{{}^tZ_J}\la|\bcb^2_{J,y}|\ra\Bpq 
\{V_{y,w,C,j}\ot{}'V_{t,w,j'}^{C',C}\ot C'_t\la j+j'\ra;\\&
w\in{}^{{}^\d J}W,C\in\uCS(\hL^{w1}),C'\in\uCS(\hL^{t1}),j\in\ZZ,j'\in\NN\}\tag e\endalign$$
in $\cd_m({}^tZ_J)$.

\subhead 3.6\endsubhead
Let $J\sub I$. For any $w\in{}^{{}^\d J}W,C\in\uCS(\hL^{w1})$ and any $K\in\cd_m(Z_J)$ we set
$$(C_w^\sha:K)=\sum_{j,h\in\ZZ}(-1)^j(\text{ multiplicity of $C_w^\sha$ in }{}^pH^j(K)_h)(-v)^h\in\ca.$$
Here, for any mixed perverse sheaf $R$, $R_h$ denotes the subquotient of $R$ of pure weight $h$. We set 
$$\c(K)=\sum_{w\in{}^{{}^\d J}W,C\in\uCS(\hL^{w1})}(C_w^\sha:K)C_w^\sha\in\fK_J.\tag a$$
For any $w\in{}^{{}^\d J}W,C\in\uCS(\hL^{w1})$ and any $K'\in\cd_m({}^wZ_J)$ we set
$$(C_w:K')=\sum_{j,h\in\ZZ}(-1)^j(\text{ multiplicity of $C_w$ in } {}^pH^j(K')_h)(-v)^h\in\ca.\tag b$$
For $K\in\cd_m(Z_J)$ we set 
$$\c_w(K)=\sum_{C\in\uCS(\hL^{w1})}(C_w:K|_{{}^wZ_J})C_w\in{}^w\fK_J.$$
From 2.10(b),(c) we deduce for $y\in W$:
$$\c_w(\bKK^y_J\la|\bcb^2_{J,y}|\ra)v^{|\bcb^2_{J,y}|}=
       \c_w(\bKK^y_J)=\sum_{y'\in W;y'\le y}P_{y',y}(v^2)\c_w(\KK^{y'}_J).\tag c$$
(We use that $\cv_{y',y}$ in 2.10(b) are pure of weight $0$, see \cite{\KLL}.)

For $t\in N_{J,\d}$, $w\in N_{J,\d}$, $C\in\uCS(\hL^{w1})$, we have
$$\c_t(C^\sha_w)=\sum_{C'\in\uCS(\hL^{t1})}\Ph_{t,w}^{C',C}C'_t\tag d$$
where 
$$\Ph_{t,w}^{C',C}=\sum_{j'\in\NN}\dim{}'V_{t,w,j'}^{C',C}v^{-j'}\in\NN[v\i]$$
specializes to the polynomial $P_{y,w}$ of \cite{\KL} (in the case where $J=\emp$, $\d=1$).
Note that $\Ph_{t,w}^{C',C}$ is $1$ if $t=w,C=C'$ and is in $v\i\NN[v\i]$ otherwise.

For $t\in N_{J,\d}$, $y\in W$ we set
$$[y]_{(t)}=v^{-|{}^tZ_J|}\c_t(\bKK^y_J).\tag e$$

\subhead 3.7\endsubhead
Let $J\sub I$, $z\in N_{J,\d}$, $y\in W$. 
For $u\in W_{{}^\d J}$  we define $\bKK^{u;\hL^{z1}}$ like $\bKK^y_J$ by replacing $\hG,G^1,J,y$ by
$\hL^z,\hL^{z1},{}^\d J,u$. We define 
$\c(\KK^{u;\hL^{z1}})\in\fK(\hL^{z1})$, $\c(\bKK^{u;\hL^{z1}})\in\fK(\hL^{z1})$, 
like $\c(\KK^y_J)$ and $\c(\bKK^y_J)$ by replacing $\hG,G^1,J,y$ by $\hL^z,\hL^{z1},{}^\d J,u$. 

From the proof of 2.5 we obtain:

If $y\n z\i W_{{}^\d J}$ then $\c_z(\KK^y_J)=0$. If $y=z\i u$ with $u\in W_{{}^\d J}$ then
$$\c_z(\KK^y_J)=v^{|{}^zZ_J|}\sum_{C\in\uCS(\hL^{z1})}f_{u,C}C_z$$
where $f_{u,C}\in\ca$ are given by  
$\c(\KK^{u;\hL^{z1}})=v^{|L^z|}\sum_{C\in\uCS(\hL^{z1})}f_{u,C}C$. Thus we have  
$$v^{-|{}^zZ_J|}\c_z(\KK^y_J)=v^{-|L^z|}(\c(\KK^{u;\hL^{z1}}))_z.\tag a$$
Let $z,t\in N_J$. Using (a), 3.6(c), we see that for $u\in W_{{}^\d J}$ we have
$$\c_t(\bKK^{z\i u}_J)=\sum_{u'\in W_{{}^\d J};t\i u'\le z\i u}P_{t\i u',z\i u}(v^2)v^{|{}^tZ_J|-|L^t|}
\c(\KK^{u';\hL^{t1}})_t.\tag b$$ 
Using 3.6(c) for $\hL^{z1}$ instead of $G^1$ we see that for $u'\in W_{{}^\d J}$ we have
$$\c(\bKK^{u';\hL^{t1}})=\sum_{u''\in W_{{}^\d J};u''\le u'}P_{u'',u'}(v^2)\c(\KK^{u'';\hL^{t1}})$$   
hence
$$\c(\KK^{u';\hL^{t1}})=\sum_{u''\in W_{{}^\d J};u''\le u'}P'_{u'',u''}(v^2)\c(\bKK^{u'';\hL^{t1}})$$
where $(P'_{u'',u'}(v^2))$ is the matrix inverse to $(P_{u'',u'}(v^2))$ (with $u',u''$ in $W_{{}^\d J}$).
Introducing this in (b) we obtain
$$[z\i u]_{(t)}=\sum_{u''\in W_{{}^\d J}}
(\sum_{u'\in W_J;u''\le u',t\i u'\le z\i u}P'_{u'',u'}(v^2)P_{t\i u',z\i u}(v^2))
\c(\bKK^{u'';\hL^{t1}})_t.\tag c$$

\subhead 3.8\endsubhead
Let $H$ be the Hecke algebra of $W$. As an $\ca$-module, $H$ has a basis $\{T_w;w\in W\}$. The 
multiplication satisfies $(T_s+1)(T_s-v^2)=0$ if $s\in I$, $T_wT_{w'}=T_{ww'}$ if $l(ww')=l(w)+l(w')$.
For any $J\sub I$ we define an $\ca$-linear map $\mu_J:H@>>>H$ by 
$$\mu_J(T_y)=T_{\d\i(y_*)}T_{y^*}\tag a$$ 
for any $y\in W$ where $y^*\in W^{{}^\d J},y_*\in W_{{}^\d J}$ are uniquely determined by $y=y^*y_*$.

Now let $w\in{}^{{}^\d J}W$ and let $(J_n,w_n)_{n\ge0}=\k(w)$. Let $n_0\ge0$ be the smallest integer such that
$J_{n_0}=J_{n_0+1}=\do$.
Let $H_w$ be the Hecke algebra defined in terms of $W_{{}^\d(J^w_\iy)}$ in the same way as $H$ was defined in 
terms of $W$. The automorphism $\t_w:W_{{}^\d(J^w_\iy)}@>>>W_{{}^\d(J^w_\iy)}$ induces an algebra automorphism 
$H_w@>>>H_w$ ($T_y\m T_{\t_w(y)}$ for $y\in W_{{}^\d(J^w_\iy)}$) denoted again by $\t_w$.

We define a $\ca$-linear function $e_w:H@>>>H_w$ by $e_w(T_y)=T_{y_1}$ if $y=w\i y_1$, 
$y_1\in W_{{}^\d(J^w_\iy)}$ and $e_w(T_y)=0$ for all other $y\in W$. For $n\ge n_0$ we define an $\ca$-linear 
function $E_{w,J,n}:H@>>>H_w$ by 
$$E_{w,J,n}(T_y)=e_w(\mu_{J_{n-1}}\do\mu_{J_2}\mu_{J_1}\mu_{J_0}(T_y)),\qua(y\in W).\tag b$$
(If $n=n_0=0$ we have $E_{w,J,n}(T_y)=e_w(T_y)$.) We show that for $y\in W$ and $n\ge n_0$ we have
$$E_{w,J,n+1}(T_y)=\t_w(E_{w,J,n}(T_y)).\tag c$$
It is enough to show that for $y\in W$ we have
$$e_w(\mu_{J^w_\iy}(T_y))=\t_w(e_w(T_y)).\tag d$$
We have $y=y'y''$ where $y'\in W^{{}^\d(J^w_\iy)},y''\in W_{{}^\d(J^w_\iy)}$ and
$$e_w(\mu_{J^w_\iy}(T_y))=e_w(T_{\d\i(y'')}T_{y'}),\qua e_w(T_y)=\d_{y',w\i}T_{y''}.$$
Now $T_{\d\i(y'')}T_{y'}$ is an $\ca$-linear combination of terms $T_{\ty y'}$ with  $\ty\in W_{J^w_\iy}$.
If $e_w(T_{\ty y'})\ne0$ then $\ty y'\in w\i W_{{}^\d(J^w_\iy)}=W_{J^w_\iy}w\i$ hence $y'=w\i$. Thus, if 
$y'\ne w\i$, both sides of (d) are zero. Now assume that $y'=w\i$. We have
$$e_w(T_{\d\i(y'')}T_{y'})=e_w(T_{\d\i(y'')}T_{w\i})=e_w(T_{w\i}T_{\t_w(y'')})=T_{\t_w(y'')}=\t_w(e_w(T_y)).
$$
This proves (d) hence also (c).

\proclaim{Proposition 3.9}Let $J\sub I$ and let $w\in{}^{{}^\d J}W$; let $n\ge n_0$ ($n_0$ as in 3.8). 
Recall the isomorphism $\fK(\hL^{w1})@>\si>>{}^w\fK_J$, $\x\m\x_w$ in 2.15. For $y\in W$ we have
$$\c_w(\KK^y_J)=v^{\ph(w,J)}\sum_{y'\in W_{{}^\d(J^w_\iy)}}f_{y',n}\c(\KK^{y',\hL^{w1}})_w\tag a$$   
where $\ph(w,J)\in\ZZ$ is explicitly computable and $f_{y',n}\in\ca$ are given (see 3.8) by
$$E_{w,J,n}(T_y)=\sum_{y'\in W_{{}^\d(J^w_\iy)}}f_{y',n}T_{y'}\in H_w.$$
\endproclaim
We argue by induction on $\sha J$. Let $z=\min(W_{{}^\d J}wW_J)$. 
Assume first that $z\in N_{J,\d}$ so that $z=w$ and $n_0=0$. If $y\n z\i W_{{}^\d J}$ then, as in the proof 
of 2.8 we have $\KK^y|_{{}^wZ_J}=0$ hence $\c_w(\KK^y_J)=0$. From the definition, in this case we have 
$E_{z,J,0}(T_y)=0$. Using 3.8(c) we deduce that $E_{z,J,n}(T_y)=0$ for any $n\ge0$ hence (a) holds in this 
case. Thus we can assume that $y=z\i y'$ with $y'\in W_{{}^\d J}$. From the definition, in this case we have 
$E_{z,J,0}(T_y)=T_{y'}$. Using 3.8(c) we deduce that for $n\ge0$ we have 
$E_{z,J,n}(T_y)=\t_w^n(T_{y'})=T_{\t_w^n(y')}$. Using 3.7(a) we have 
$\c_z(\KK^{z\i y'}_J)=v^{\ph(z,J)}\c(\KK^{y',\hL^{z1}})_z$ where $\ph(z,J)=|{}^zZ_J|-|L^z|$.
 It remains use that
$\c_z(\KK^{z\i y'}_J)=\c_z(\KK^{z\i\t_z^n(y')}_J)$ which follows from 2.7(c).

Next we assume that $z\n N_{J,\d}$. Then, setting $J_1=J\cap\d\i(\Ad(z)J)$, we have $\sha J_1<\sha J$. We
can write uniquely $y=y^*y_*$ as in 2.7. From 2.7(b) we see that
$$\c_w(\KK^y_J)=v^{\ph'(w,J)}\sum_{y'\in W_J}g_{y'}\c_w(\KK^{y'y^*}_{J_1})$$
where $T_{\d\i(y_*)}T_{y^*}=\sum_{y'\in W_J}g_{y'}T_{y'y^*}$ in $H$ (with $g_{y'}\in\ca$) and 
$\ph'(w,J)\in\ZZ$. By the induction hypothesis we have
$$\c_w(\KK^{y'y^*}_{J_1})=v^{\ph(w,J_1)}
\sum_{\ty'\in W_{{}^\d(J^w_\iy)}}\tf_{\ty',n-1}\c(\KK^{\ty',\hL^{w1}})_w,$$
where
$$E_{w,J_1,n-1}(T_{y'y^*})=\sum_{\ty'\in W_{{}^\d(J^w_\iy)}}\tf_{\ty',n-1}T_{\ty'}\in H_w$$
and $\ph(w,J_1)\in\ZZ$. It follows that
$$\c_w(\KK^y_J)=v^{\ph'(w,J)+\ph(w,J_1)}\sum_{y'\in W_J,\ty'\in W_{{}^\d(J^w_\iy)}}g_{y'}\tf_{\ty',n-1}
\c(\KK^{\ty',\hL^{w1}})_w.$$
We have 
$$\align&\sum_{y'\in W_J,\ty'\in W_{{}^\d(J^w_\iy)}}g_{y'}\tf_{\ty',n-1}T_{\ty'}=
\sum_{y'\in W_J}g_{y'}E_{w,J_1,n-1}(T_{y'y^*})\\&
=\sum_{y'\in W_J}g_{y'}e_w(\mu_{J_{n-1}}\do\mu_{J_2}\mu_{J_1}(T_{y'y^*}))
=e_w(\mu_{J_{n-1}}\do\mu_{J_2}\mu_{J_1}(T_{\d\i(y_*)}T_{y^*}))\\&
=e_w(\mu_{J_{n-1}}\do\mu_{J_2}\mu_{J_1}\mu_J(T_y))=E_{w,J,n}(T_y).\endalign$$
Thus (a) holds with $\ph(w,J)=\ph'(w,J)+\ph(w,J_1)$.

\proclaim{Proposition 3.10} Let $J\sub I$ and let $\o_J\in W_J$ be as in 3.2. Let $w\in{}^{{}^\d J}W$.

(a) If $z\in N_{J,\d}$, then $\bKK^{w\i\d(\o_J)}_J|_{{}^zZ_J}$ is a direct sum of complexes of the form 
$\bbq\la j\ra$ with $j\in\ZZ$.

(b) Let $C=\bbq\la|L^w|\ra\in CS(\hL^{w1})$. Then for some $j\in\ZZ$, $C^\sha_w\la j\ra$ is a direct summand of $\bKK^{w\i\d(\o_J)}_J$.

(c) If $C$ is as in (b) and $z$ is as in (a) then $C^\sha_w|_{{}^zZ_J}$ is a direct sum of complexes of the 
form $\bbq\la j\ra$ with $j\in\ZZ$.
\endproclaim
We prove (a). Using 3.4(a) and 3.6(c) we see that it is enough to show that
$$\sum_{y'\in W;y'\le w\i\d(\o_J)}P_{y',w\i\o_J}(v^2)\c_z(\KK^{y'}_J)\in\ca \bbq\la|{}^zZ_J|\ra.$$
Since $w\i\in W^{{}^\d J}$, the last sum is equal to
$$\sum_{y'\in W^{{}^\d J};y'\le w\i\d(\o_J)}P_{y',w\i\o_J}(v^2)
\sum_{u\in W_{{}^\d(J)}}\c_z(\KK^{y'u}_J).$$
(We use a standard property of the polynomials $P_{y',y}$.) Hence it is enough to show that
$$\sum_{u\in W_{{}^\d(J)}}\c_z(\KK^{y'u}_J)\in\ca \bbq\la|{}^zZ_J|\ra$$
for any $y'\in W^{{}^\d J}$. By arguments in 3.7, the left hand side is zero unless $y'=z\i$ in which case it
equals
$$v^{|{}^zZ_J|-|L^z|}\sum_{u\in W_{{}^\d(J)}}(\c(\KK^{u;\hL^{z1}}))_z=
v^{|{}^zZ_J|-|L^z|}(\c(\bKK^{\d(\o_J);\hL^{z1}}))_z.$$
It remains to use that
$$\c(\bKK^{\d(\o_J);\hL^{z1}})\in\ca \bbq\la|\hL^{z1}|\ra.$$

We omit the proof of (b). Now (c) follows immediately from (a) and (b).

\subhead 3.11\endsubhead
In this subsection we assume that $\hG=G$ is simple modulo its centre hence $G^1=G$ and $\d=1$. We fix $J\sub I$ and 
we write $N_J$ instead of $N_{J,\d}=\{w\in W;wJw\i=J\}$, a subgroup of $W$ with unit element $e$. We assume 
that $J\ne I$ and that we are given a cuspidal object $A_1$ of $\uCS(L^e)$ (with $L^e$ as in 3.1). It 
follows that $L^e$ modulo its centre is simple or $\{1\}$. We will attach to each $w\in N_J$ a (not 
identically zero) map $\io_w:\uCS(\hL^{w1})@>>>\ZZ$ (said to be a {\it cuspidal map}) well defined up to
multiplication by $\pm1$. Let $F_w:L^e@>>>L^e$ be the Frobenius map for an $\FF_q$-rational structure on 
$L^e$ whose action on the Weyl group $W_J$ of $L^e$ (or of $L^w$) is induced by the conjugation action of 
the connected component $\hL^{w1}$ of 
$\hL^w$ on $L^w$ (see 2.15). As in \cite{\URCC} we identify $\uCS(\hL^{w1})$ with $\ucs(L^e,F_w)$, a set of 
representatives for the isomorphism classes of unipotent representations of the finite reductive group 
$(L^e)^{F_w}$. Then $\io_w$ becomes a map $\ucs(L^e,F_w)@>>>\ZZ$. Now $A_1$ gives rise to a cuspidal object
$A_w$ of $\ucs(L^e,F_e)$ and as in \cite{\URCC} this corresponds to a unipotent cuspidal representation
$\r_w$ of $(L^w)^{F_e}$. According to \cite{\ORA, 4.23}, $\r_w$ has an
associated two-sided cell $\boc$ of $W_J$ and it corresponds to a pair $(x,r)$ where $x$ is an element of
a certain finite group $\G$ attached to $\boc$ and $r$ is an irreducible representation of $Z_\G(x)$.
Moreover, $\boc$ gives rise to a subset $\uCS(\hL^{w1})_\boc$ of $\uCS(\hL^{w1})$ or equivalently
to a subset $\ucs(L^e,F_w)_\boc$ of $\ucs(L^e,F_w)$ in natural bijection \cite{\ORA, 4.23} with the
set $\bar M(\G\sub\ti\G)$ defined as follows: $\G\sub\ti\G$ is a certain imbedding of $\G$
as a normal subgroup into a finite group $\ti\G$ such that $\ti\G/\G$ is cyclic of order $\le3$ with a given
$\G$-coset $\G^1$ whose image in $\ti\G/\G$ generates $\ti\G/\G$; 
$\bar M(\G\sub\ti\G)$ consists of all pairs $(y,s)$ where $y\in\G^1$ is defined up to 
conjugation in $\ti\G$ and $s$ is an irreducible representation of $Z_\G(y)$, the centralizer of $y$ in $\G$,
up to isomorphism. Our function $\io_w$ is required to be zero on the complement of $\ucs(L^e,F_w)_\boc$
hence it can be viewed as a function $\io_w:\bar M(\G\sub\ti\G)@>>>\ZZ$.
Let $M(\ti\G)$ be the set consisting of all pairs $(y',s')$ where $y'\in\ti\G$ is defined up to 
conjugation in $\ti\G$ and $s'$ is an irreducible representation of $Z_{\ti\G}(y')$, the centralizer of $y'$ 
in $\ti\G$, up to isomorphism. We can find an irreducible representation $r'$ of $Z_{\ti\G}(x)$ whose 
restriction to $Z_\G(x)$ is $r$.

Let $\{,\}:M(\ti\G)\T M(\ti\G)@>>>\bbq$ be the pairing \cite{\ORA, 4.14.3}.
We define $j_{r'}:\bar M(\G\sub\ti\G)@>>>\ZZ$ by $j_{r'}(y,s)=\{(x,r'),(y,s')\}$ where
$s'$ is an irreducible representation of $Z_{\ti\G}(y)$ whose restriction to $Z_\G(y)$ is $s$. Since
$x\in\G$, $j_{r'}(y,s)$ is independent of the choice of $s'$. We can choose $r'$ so that $j_{r'}(y,s)$ 
takes values in $\QQ$. We define $\io_w:\bar M(\G\sub\ti\G)@>>>\ZZ$ as $cj_{r'}$ where $c$ is a rational 
number $>0$ such that $\io_w$ takes values in $\ZZ$ and there is no integer $d>1$ dividing all values of
$\io_w$. In the case where $\ti\G/\G$ has order $1$ or $3$, $\io_w$ is unique.
In the case where $\ti\G/\G$ has order $2$, $\io_w$ is unique up to multiplication by $\pm1$. We state some
conjectures.

{\it Conjecture 1. For any $t\le z$ in $N_J$ there is a (necessarily unique) $X_{t,z}\in\ZZ[v\i]$ such that 
$$\sum_{C\in\uCS(\hL^{z1})_\boc}\io_z(C)\Ph_{t,z}^{C',C}=X_{t,z}\io_t(C')$$
for any $C'\in\uCS(\hL^{t1})_\boc$ where $\Ph_{t,z}^{C',C}\in\NN[v\i]$ are as in 3.6.}

An equivalent statement is that in ${}^t\fK_J$ we have
$$\sum_{C\in\uCS(\hL^{z1})_\boc}\io_z(C)\c_t(C^\sha_z)=X_{t,z}\sum_{C'\in\uCS(\hL^{t1})_\boc}
\io_t(C')C'_t$$
modulo $\sum_{C'\in\uCS(\hL^{t1})-\uCS(\hL^{t1})_\boc}\ca C'_t$.

{\it Conjecture 2. For any $t\le z$ in $N_J$, the matrix
$(\Ph_{t,z}^{C',C})_{(C',C)\in\uCS(\hL^{t1})_\boc\T\uCS(\hL^{z1})_\boc}$ is square and invertible.}

For any $i\in I-J$ we have $i\in N_J$. It follows that $\o_{J\cup i}\o_J=\o_J\o_{J\cup i}\in N_J$ hence
$\t_i=\o_{J\cup i}\o_J=\o_J\o_{J\cup i}$ has square $1$. It is known that $N_J$ together with
$\{\t_i;i\in I-J\}$ is a Coxeter group. Let $a:W@>>>\NN$ be the $a$-function of $W,l$, see \cite{\HEC, 13.6}.
For $i\in I-J$ we set $c_i=a(x\t_i)-a(x')$ where $x,x'\in\boc$; this independent of the choice of $x,x'$
by \cite{\HEC, 9.13, 4.12(P11)}. There is a unique weight function $\cl:N_J@>>>\NN$ such that
$\cl(\t_i)=c_i$ for all $i\in I-J$. Hence the Hecke algebra $\ch$ associated to $N_J,\cl$ and the elements
$p_{t,z}$ (for $t,z$ in $N_J$) are well defined as in 0.1.

{\it Conjecture 3. For any $t\le z$ in $N_J$ we have $X_{t,z}=e_{t,z}p_{t,z}$ where $e_{t,w}=\pm1$.}

\head 4. An example \endhead
\subhead 4.1\endsubhead
In this section we assume that we are in the setup of 3.1 and that $\hG=G^1=G$ is such that $W$ is of 
type $B_4$. We have $\d=1$.
We shall denote the elements of $I$ as $s_i$ ($i=1,2,3,4$) where the notation is chosen so 
that $(s_1s_2)^4=(s_2s_3)^3=(s_3s_4)^3=1$ and $s_is_j=s_js_i$ if $i-j\in\{\pm2,\pm3\}$. An element $w\in W$ 
with reduced expression $s_{i_1}s_{i_2}\do s_{i_m}$ will be denoted as $i_1i_2\do i_m$. In particular we 
write $i$ instead of $s_i$; the unit element of $W$ is denoted by $\emp$. We set $J=\{1,2\}\sub I$. 
The elements of ${}^JW^J$ are 
$$\align&\emp,3,4,34,43,343,3243,32123,321234,321243,432123,\\&3212343,3432123,4321234,34321234,32123432123,
321234321234.\endalign$$
Now $N_J:=N_{J,1}$ is the subgroup of $W$ consisting of the elements
$$\align&\emp,e=4,f=32123,fe=321234,ef=432123,\\&efe=4321234,fef=32123432123,efef=fefe=321234321234.
\endalign$$
It is a Coxeter group (of order $8$) with generators $e$, $f$ which satisfy $(ef)^4=1$. We define a weight 
function $\cl:N_J@>>>\NN$ by
$$\emp\m0,e\m1,f\m3,fe\m4,ef\m4,efe\m5,fef\m7,efef\m8$$
and a homomorphism $\e:N_J\m\{\pm1\}$ by $\e(e)=1$, $\e(f)=-1$.
Note that $\cl$ coincides with the weight function defined in 3.11 in terms of $W, W_J$ and the two-sided
cell $\boc=\{1,2,12,21,121,212\}$ of $W_J$.

If $z\in N_J$ then the Weyl group of $L^z$ is $W_J=\{\emp,1,2,12,21,121,212,1212\}$. In our case we have 
$L^z=\hL^{z1}$. Also $L^z$ is independent of $z$ (in $N_J$) up to an inner automorphism; hence we can use 
the notation $L$ instead of $L^z$. For $u\in W_J$ we set $[u]=v^{-|L|}\c(\bKK^{u;L})$. 

\subhead 4.2\endsubhead
The objects of $\uCS(L)$ can be denoted by $1,\r,\s,\s',\th,S$. Here $1,\r,\s,\s',S$ are perverse sheaves on
$L$ with support equal to $L$ which are generically local systems (up to shift) of rank $1,2,1,1,1$; $\th$ 
is a cuspidal character sheaf on $L$. They can be characterized by the equalities

$[\emp]=1+2\r+\s+\s'+S$;

$[1]=(1+v^2)(1+\r+\s)$;

$[2]=(1+v^2)(1+\r+\s')$;

$[12]=v^2(\r+\s+\s'+\th)+(1+v^2)^21$;

$[21]=v^2(\r+\s+\s'+\th)+(1+v^2)^21$;

$[121]=(v^2+v^4)(\s'+\th)+(1+v^2)(1+v^4)1$;

$[212]=(v^2+v^4)(\s+\th)+(1+v^2)(1+v^4)1$;

$[1212]=(1+v^2)^2(1+v^4)1$.
\nl
We have $1=\bbq\la|L|\ra$.

\subhead 4.3\endsubhead
Let $z,t\in N_J$, $u\in W_J$.
In our case the explicit values of $P_{t\i u',z\i u}(v^2)$ in 3.7(c) can be found in the tables of 
\cite{\GOR}; moreover in 3.7(c) we have $P'_{u',u''}(v^2)=(-1)^{l(u')+l(u'')}$. Hence 
the coefficients of $[z\i u]_{(t)}$ in 3.7(c) are explicitly known.
In subsections 4.4-4.10 we give for any $z,t$ in $N_J$ the explicit values of $[z\i u]_{(t)}$, with 
$u\in W_J-\{\emp,1212\}$, as an $\ca$-linear combination of elements $[u'']_t$ with $u''\in W_J$. For 
$\x,\x'$ in ${}^t\fK_J$ we write $\x\si\x'$ if $\x-\x'\in\ca1_t$.

\subhead 4.4\endsubhead
Assume that $(t,z)$ satisfies either $t=z$ or that it is one of 
$$\align&(\emp,4), (32123,432123), (32123,321234),(32123,4321234),\\&(432123,4321234),(321234,4321234),
(32123432123,321234321234)\endalign$$
that is,
$$(\emp,e),(f,ef),(f,fe),(f,efe),(ef,efe),(fe,efe),(fef,efef).$$
Note that $l(z)-l(t)=\cl(z)-\cl(t)$, $\e(z)\e(t)=1$. From 3.7(c) we have

$[z\i(121)]_{(t)}=[121]_t$,

$[z\i(212)]_{(t)}=[212]_t$,

$[z\i(12)]_{(t)}=[12]_t$,

$[z\i(21)]_{(t)}=[21]_t$,

$[z\i(2)]_{(t)}=[2]_t$,

$[z\i(1)]_{(t)}=[1]_t$.

Hence, using 4.2, we have

$[z\i(121)]_{(t)}\si(v^2+v^4)(\th_t+\s'_t)$,

$[z\i(212)]_{(t)}\si(v^2+v^4)(\th_t+\s_t)$,

$[z\i(12)]_{(t)}\si v^2(\r_t+\s_t+\s'_t+\th_t)$,

$[z\i(21)]_{(t)}\si v^2(\r_t+\s_t+\s'_t+\th_t)$,

$[z\i(2)]_{(t)}\si (1+v^2)(\r_t+\s'_t)$,

$[z\i(1)]_{(t)}\si (1+v^2)(\r_t+\s_t)$.

\subhead 4.5\endsubhead
Assume that $(t,z)$ is one of 
$$\align&(\emp,32123),(\emp,432123), (4,432123),(\emp,321234),(4,321234),\\&
(32123,321234321234),(432123,321234321234),(321234,321234321234),\\&(4321234,321234321234),
(432123,32123432123), (321234,32123432123)\endalign$$
that is, one of
$$\align&(\emp,f),(\emp,ef),(e,ef), (\emp,fe),(e,fe),\\& 
(f,efef), (ef,efef), (fe,efef), (efe,efef), (ef, fef), (fe,fef).\endalign$$
Note that $l(z)-l(t)=\cl(z)-\cl(t)+2$, $\e(z)\e(t)=-1$. From 3.7(c) we have

$[z\i(121)]_{(t)}=v^2[121]_t+[1212]_t$,

$[z\i(212)]_{(t)}=v^4[2]_t+[1212]_t$,

$[z\i(12)]_{(t)}=v^2[12]_t+[1212]_t$,

$[z\i(21)]_{(t)}=v^2[21]_t+[1212]_t$,

$[z\i(2)]_{(t)}=[212]_t$,

$[z\i(1)]_{(t)}=v^2[1]_t+[1212]_t$.

Hence, using 4.2, we have

$[z\i(121)]_{(t)}\si(v^4+v^6)(\th_t+\s'_t)$,

$[z\i(212)]_{(t)}\si(v^4+v^6)(\r_t+\s'_t)$,

$[z\i(12)]_{(t)}\si v^4(\r_t+\s_t+\s'_t+\th_t)$,

$[z\i(21)]_{(t)}\si v^4(\r_t+\s_t+\s'_t+\th_t)$,

$[z\i(2)]_{(t)}\si(v^2+v^4)(\th_t+\s_t)$,

$[z\i(1)]_{(t)}\si(v^2+v^4)(\r_t+\s_t)$.

\subhead 4.6\endsubhead
Assume that $(t,z)$ is one of $(\emp,4321234), (4,4321234)$ that is, one of \lb $(\emp,efe),(e,efe)$. 

Note that $l(z)-l(t)=\cl(z)-\cl(t)+2=-1$, $\e(z)\e(t)=-1$. From 3.7(c) we have

$[z\i(121)]_{(t)}=(v^2+v^4)[121]_t+(1+v^2)[1212]_t$,

$[z\i(212)]_{(t)}=(v^4+v^6)[2]_t+(1+v^2)[1212]_t$,

$[z\i(12)]_{(t)}=(v^2+v^4)[12]_t+(1+v^2)[1212]_t$,

$[z\i(21)]_{(t)}=(v^2+v^4)[21]_t+(1+v^2)[1212]_t$,

$[z\i(2)]_{(t)}=(1+v^2)[212]_t$,

$[z\i(1)]_{(t)}=(v^2+v^4)[1]_t+(1+v^2)[1212]_t$.

Hence, using 4.2, we have

$[z\i(121)]_{(t)}\si(v^2+v^4)^2(\th_t+\s'_t)$,

$[z\i(212)]_{(t)}\si(v^4+v^6)(1+v^2)(\r_t+\s'_t)$,

$[z\i(12)]_{(t)}\si(v^4+v^6)(\r_t+\s_t+\s'_t+\th_t)$,

$[z\i(21)]_{(t)}\si(v^4+v^6)(\r_t+\s_t+\s'_t+\th_t)$,

$[z\i(2)]_{(t)}\si(1+v^2)(v^2+v^4)(\th_t+\s_t)$,

$[z\i(1)]_{(t)}\si(1+v^2)(v^2+v^4)(\r_t+\s_t)$.

\subhead 4.7\endsubhead
Assume that $(t,z)$ is one of 
$$(\emp,321234321234), (4,321234321234),(4,32123432123)$$ 
that is, one of $(\emp,efef),(e,efef),(e,fef)$. 

Note that $l(z)-l(t)=\cl(z)-\cl(t)+4$, $\e(z)\e(t)=1$. From 3.7(c) we have

$[z\i(121)]_{(t)}=v^4[121]_t$,

$[z\i(212)]_{(t)}=v^4[212]_t$,

$[z\i(12)]_{(t)}=v^4[12]_t$,

$[z\i(21)]_{(t)}=v^4[21]_t$,

$[z\i(2)]_{(t)}=v^4[2]_t$,

$[z\i(1)]_{(t)}=v^4[1]_t$.

Hence, using 4.2, we have

$[z\i(121)]_{(t)}\si(v^6+v^8)(\th_t+\s'_t)$,

$[z\i(212)]_{(t)}\si(v^6+v^8)(\th_t+\s_t)$,

$[z\i(12)]_{(t)}\si v^6(\r_t+\s_t+\s'_t+\th_t)$,

$[z\i(21)]_{(t)}\si v^6(\r_t+\s_t+\s'_t+\th_t)$,

$[z\i(2)]_{(t)}\si(v^4+v^6)(\r_t+\s'_t)$,

$[z\i(1)]_{(t)}\si(v^4+v^6)(\r_t+\s_t)$.

\subhead 4.8\endsubhead
Assume that $(t,z)=(\emp,32123432123)$ that is, $(\emp,fef)$. 

Note that $l(z)-l(t)=\cl(z)-\cl(t)+4$, $\e(z)\e(t)=1$. From 3.7(c) we have

$[z\i(121)]_{(t)}=v^8[1]_t+v^4[121]_t$,

$[z\i(212)]_{(t)}=(v^4+v^6)[212]_t$,

$[z\i(12)]_{(t)}=(v^4+v^6)[12]_t$,

$[z\i(21)]_{(t)}=(v^4+v^6)[21]_t$,

$[z\i(2)]_{(t)}=(v^4+v^6)[2]_t$,

$[z\i(1)]_{(t)}=v^4[121]_t+v^4[1]_t$.

Hence, using 4.2, we have

$[z\i(121)]_{(t)}\si(v^8+v^{10})(\r_t+\s_t)+(v^6+v^8)(\th_t+\s'_t)$,

$[z\i(212)]_{(t)}\si(v^4+v^6)(v^2+v^4)(\th_t+\s_t)$,

$[z\i(12)]_{(t)}\si(v^6+v^8)(\r_t+\s_t+\s'_t+\th_t)$,

$[z\i(21)]_{(t)}\si(v^6+v^8)(\r_t+\s_t+\s'_t+\th_t)$,

$[z\i(2)]_{(t)}\si(v^4+v^6)(1+v^2)(\r_t+\s'_t)$,

$[z\i(1)]_{(t)}\si(v^6+v^8)(\th_t+\s'_t)+(v^4+v^6)(\r_t+\s_t)$.

\subhead 4.9\endsubhead
Assume that $(t,z)=(32123,32123432123)$ that is, $(f,fef)$.

Note that $l(z)-l(t)=\cl(z)-\cl(t)+2$, $\e(z)\e(t)=-1$. From 3.7(c) we have

$[z\i(121)]_{(t)}=v^6[1]_t+v^2[121]_t$,

$[z\i(212)]_{(t)}=(v^4+v^6)[2]_t$,

$[z\i(12)]_{(t)}=(v^2+v^4)[12]_t$,

$[z\i(21)]_{(t)}=(v^2+v^4)[21]_t$,

$[z\i(2)]_{(t)}=(1+v^2)[212]_t$,

$[z\i(1)]_{(t)}=v^2[121]_t+v^2[1]_t$.

Hence, using 4.2, we have

$[z\i(121)]_{(t)}\si(v^6+v^8)(\r_t+\s_t)+(v^4+v^6)(\th_t+\s'_t)$,

$[z\i(212)]_{(t)}\si(v^4+v^6)(1+v^2)(\r_t+\s'_t)$,

$[z\i(12)]_{(t)}\si(v^4+v^6)(\r_t+\s_t+\s'_t+\th_t)$,

$[z\i(21)]_{(t)}\si(v^4+v^6)(\r_t+\s_t+\s'_t+\th_t)$,

$[z\i(2)]_{(t)}\si(1+v^2)(v^2+v^4)(\th_t+\s_t)$,

$[z\i(1)]_{(t)}\si(v^4+v^6)(\th_t+\s'_t)+(v^2+v^4)(\r_t+\s_t)$.

\subhead 4.10\endsubhead
If $(t,z)$ in $N_J$ is not as in 4.4-4.9 then we have $[z\i u]_{(t)}=0$ for any $u\in W_J$.

\subhead 4.11\endsubhead
Let $z,t\in N_J$, $u\in W_J$. We set
$$[z\i u]'_{(t)}=\c_t(\bKK^{z\i u}_J\la|\bcb^2_{J,z\i u}|\ra).$$
Using 2.10(d) we have
$$[z\i u]'_{(t)}=v^{-|\bcb^2_{J,z\i u}|+|{}^tZ_J|}[z\i u]_{(t)}=v^{-l(z)+l(t)-l(u)}[z\i u]_{(t)}.\tag a$$
Let $C\in\uCS(L)-\{1\}$, let $t,z\in N_J$ and let $u\in W_J-\{\emp,1212\}$. From 4.4-4.10 we see that the 
following result holds:

(b) {\it The coefficient of $C_t$ in $[z\i u]'_{(t)}$ is in $\NN[v\i]$. More precisely, this coefficient is:

(i) $0$ or $1+v^{-2}$ if $(t,z)$ is as in 4.4 with $t\ne z$, $u\in\{1,2,121,212\}$;

(ii) in $v\i\NN[v\i]$ if $t,z,u$ are not as in (i) but $t\ne z$.}

\subhead 4.12\endsubhead
Let $z,t\in N_J$, $u\in W_J$. Using 4.11(a) and 3.5(e) with $y=z\i u$ we deduce
$$v^{-l(z)+l(t)-l(u)}[z\i u]_{(t)}=\sum_{C'\in\uCS(L)}N^{z,t,u}_{C'}C'_t$$
where
$$N^{z,t,u}_{C'}=\sum_{w\in{}^JW,C\in\uCS(L),j\in\ZZ,j'\in\NN}
\dim V_{z\i u,w,C,j}\dim{}'V_{t,w,j'}^{C',C}v^{-j-j'}.$$
If $z=t$, in the previous sum we must have $w=z$. Note that $\dim{}'V_{z,z,j'}^{C',C}=0$ unless $C=C'$, 
$j'=0$ and we have
$$N^{z,z,u}_{C'}=\sum_{j\in\ZZ}\dim V_{z\i u,z,C',j}v^{-j}.\tag a$$
We show:

\proclaim{Proposition 4.13} Let $z,t\in N_J$ with $z\ne t$. Let $j\in\ZZ$, $u\in W_J-\{\emp,1212\}$ 
and $C'\in\uCS(L)-\{1\}$. Then $V_{z\i u,t,C',j}=0$ (notation of 3.5).
\endproclaim
We first note that if $t,z,u$ are as in 4.11(i), then from the definitions we have

(a) ${}'V_{t,z,1}^{C',C'}={}'V_{z,z,0}^{C',C'}$.
\nl
For $\x,\x'$ in $\ca$ we write $\x\ge\x'$ if $\x-\x'\in\NN[v,v\i]$. Note that $N^{z,t,u}_{C'}$ is $\ge$ 
than the corresponding sum in which $(w,C,j')$ is restricted

to $(t,C',0)$ (if $t,z,u$ are as in 4.11(ii));

to $(t,C',0)$ or to $(z,C',1)$ (if $t,z,u$ are as in 4.11(i)).
\nl
Thus
$$N^{z,t,u}_{C'}\ge\sum_{j\in\ZZ}\dim V_{z\i u,t,C',j}\dim {}'V_{t,t,0}^{C',C'}v^{-j}=
\sum_{j\in\ZZ}\dim V_{z\i u,t,C',j}v^{-j}$$
(if $t,z,u$ are as in 4.11(ii)) and
$$\align&N^{z,t,u}_{C'}\ge\sum_{j\in\ZZ}\dim V_{z\i u,t,C',j}\dim{}'V_{t,t,0}^{C',C'}v^{-j}\\&+
\sum_{j\in\ZZ}\dim V_{z\i u,z,C',j}\dim{}'V_{t,z,1}^{C',C'}v^{-j-1}
=\sum_{j\in\ZZ}\dim V_{z\i u,t,C',j}\dim{}'V_{t,t,0}^{C',C'}v^{-j}\\&+
\sum_{j\in\ZZ}\dim V_{z\i u,z,C',j}\dim {}'V_{z,z,0}^{C',C'}v^{-j-1}
=\sum_{j\in\ZZ}\dim V_{z\i u,t,C',j}v^{-j}+v\i N^{z,z,u}_{C'}\endalign$$
(if $t,z,u$ are as in 4.11(i)). (We have used (a) and 4.12(a).) If $t,z,u$ are as in 4.11(i)), we have 
$N^{z,t,u}_{C'}=v\i N^{z,z,u}_{C'}$ (see 4.4) and we deduce that 
$$0\ge\sum_{j\in\ZZ}\dim V_{z\i u,t,C',j}v^{-j};$$ 
hence $V_{z\i u,t,C',j}=0$ for all $j$. If $t,z,u$ are as in 4.11(ii)), we see using 3.5(b) that the sum 
$\sum_{j\in\ZZ}\dim V_{z\i u,t,C',j}v^{-j}$ is either zero, or

{\it for some $j$, $v^j$ and $v^{-j}$ both appear in it with $>0$ coefficient.}
\nl
In the last case it follows that $v^j$ and $v^{-j}$ both appear in $N^{z,t,u}_{C'}$ with $>0$ coefficient. 
This is not compatible with the inclusion $N^{z,t,u}_{C'}\in v\i\NN[v\i]$. Thus we must have
$\sum_{j\in\ZZ}\dim(V_{z\i u,t,C',j})v^{-j}=0$ hence $V_{z\i u,t,C',j}=0$ for all $j$. The proposition is 
proved.

\proclaim{Proposition 4.14} Let $z\in N_J$, $u\in W_J-\{\emp,1212\}$, $j\in\ZZ$. Let $t\in{}^JW$ be such 
that $\sha J^t_\iy<\sha J$. Assume that $C\in\uCS(\hL^{t1})$ is not isomorphic to $\bbq\la|L^t|\ra$. Then 
$V_{z\i u,t,C,j}=0$ (notation of 3.5(a)).
\endproclaim
The existence of $C$ guarantees that $L^t$ is not a torus. Thus $J^t_\iy$ consists of a single element $i$
(equal to $1$ or $2$) and $L^t$ has semisimple rank $1$. Hence $C$ is uniquely determined and it appears with
coefficient $1$ in $v^{-|L^t|}\c(\KK^{\emp,L^t})$ and with coefficient $-1$ in $v^{-|L^t|}\c(\KK^{i,L^t})$.
Hence the coefficient of $C_t$ in $v^{-|\bcb^2_{J,y}|}\c_t(\KK^y_J)$ is explicitly computable from 3.9(a) 
for any $y\in W$. Using now 3.6(c) (in which the polynomials $P_{y',y}(v^2)$ are explicitly known from 
\cite{\GOR}) we see that the coefficient of $C_t$ in $v^{-|\bcb^2_{J,y}|}\c_t(\bKK^y_J)$ is explicitly 
computable for any $y\in W$. In particular, the coefficient of $C_t$ in
$v^{-|\bcb^2_{J,z\i u}|}\c_t(\bKK^{z\i u}_J)$ is explicitly computable. We find that this coefficient is in 
$v\i\ZZ[v\i]$. On the other hand this coefficient is equal to
$\sum_{j\in\ZZ}\dim(V_{z\i u,t,C,j})v^{-j}$ which is invariant under the involution $v\m v\i$ of $\ca$, see
3.5(b). This forces $\dim(V_{z\i u,t,C,j})$ to be zero for any $j$. The proposition is proved.

\proclaim{Proposition 4.15}  Let $z\in N_J$, $u\in W_J-\{\emp,1212\}$.  

(a) We have
$$\bKK^{z\i u}_J\la|\bcb^2_{J,z\i u}|\ra\cong\op_{C\in\uCS(L)-\{1\},j\in\ZZ}V_{z\i u,z,C,j}\ot C_z^\sha\la 
j\ra\op K'$$
where 
$$K'\cong\op_{w\in{}^JW,j\in\ZZ}\tV_{w,j}\ot(\bbq\la|L^w|\ra)_w^\sha\la j \ra$$
and $\tV_{w,j}$ are certain $\bbq$-vector spaces.

(b) For any $t\in N_J$ we have (with $\si$ as in 4.3):
$$[z\i u]'_{(t)}\si\sum_{C\in\uCS(L)-\{1\},j\in\ZZ}\dim(V_{z\i u,z,C,j})v^{-j}\c_t(C_z^\sha).$$
\endproclaim
(a) follows from 4.13, 4.14; (b) follows from (a) using 3.10(a).   

\subhead 4.16\endsubhead
In the setup of 4.15(b), the integers $\dim(V_{z\i u,z,C,j})$ (with $C\ne1$) can be obtained from 4.4 (with 
$t=z$). Thus we can rewrite 4.15(b) as follows (recall that $z,t\in N_J$;
we set $\z=v^{-l(z)+l(t)}$ and we take $u\in\{1,2,121,212\}$):
$$\z v^{-3}[z\i121]_t\si v^{-3}(v^2+v^4)(\c_t(\th_z^\sha)+\c_t(\s'{}_z^\sha)),$$
$$\z v^{-3}[z\i212]_t\si v^{-3}(v^2+v^4)(\c_t(\th_z^\sha)+\c_t(\s_z^\sha)),$$
$$\z v^{-1}[z\i2]_t\si v\i(1+v^2)(\c_t(\r_z^\sha)+\c_t(\s'{}_z^\sha)),$$
$$\z v^{-1}[z\i1]_t\si v\i(1+v^2)(\c_t(\r_z^\sha)+\c_t(\s_z^\sha)).$$

Using now the formulas in 4.4-4.10 we deduce that the following hold.

If $(t,z)$ are as in 4.4 then
$$\z v^{-3}(v^2+v^4)(\th_t+\s'_t)\si v^{-3}(v^2+v^4)(\c_t(\th_z^\sha)+\c_t(\s'{}_z^\sha)),$$
$$\z v^{-3}(v^2+v^4)(\th_t+\s_t)\si v^{-3}(v^2+v^4)(\c_t(\th_z^\sha)+\c_t(\s_z^\sha)),$$
$$\z v^{-1}(1+v^2)(\r_t+\s'_t)\si v\i(1+v^2)(\c_t(\r_z^\sha)+\c_t(\s'{}_z^\sha)),$$
$$\z v^{-1}(1+v^2)(\r_t+\s_t)\si v\i(1+v^2)(\c_t(\r_z^\sha)+\c_t(\s_z^\sha)).$$

If $(t,z)$ are as in 4.5 then
$$\z v^{-3}(v^4+v^6)(\th_t+\s'_t)\si v^{-3}(v^2+v^4)(\c_t(\th_z^\sha)+\c_t(\s'{}_z^\sha)),$$
$$\z v^{-3}(v^4+v^6)(\r_t+\s'_t)\si v^{-3}(v^2+v^4)(\c_t(\th_z^\sha)+\c_t(\s_z^\sha)),$$
$$\z v^{-1}(v^2+v^4)(\th_t+\s_t)\si v\i(1+v^2)(\c_t(\r_z^\sha)+\c_t(\s'{}_z^\sha)),$$
$$\z v^{-1}(v^2+v^4)(\r_t+\s_t)\si v\i(1+v^2)(\c_t(\r_z^\sha)+\c_t(\s_z^\sha)).$$

If $(t,z)$ are as in 4.6 then
$$\z v^{-3}(v^2+v^4)^2(\th_t+\s'_t)\si v^{-3}(v^2+v^4)(\c_t(\th_z^\sha)+\c_t(\s'{}_z^\sha)),$$
$$\z v^{-3}(v^4+v^6)(1+v^2)(\r_t+\s'_t)\si v^{-3}(v^2+v^4)(\c_t(\th_z^\sha)+\c_t(\s_z^\sha)),$$
$$\z v^{-1}(1+v^2)(v^2+v^4)(\th_t+\s_t)\si v\i(1+v^2)(\c_t(\r_z^\sha)+\c_t(\s'{}_z^\sha)),$$
$$\z v^{-1}(1+v^2)(v^2+v^4)(\r_t+\s_t)\si v\i(1+v^2)(\c_t(\r_z^\sha)+\c_t(\s_z^\sha)).$$

If $(t,z)$ are as in 4.7 then
$$\z v^{-3}(v^6+v^8)(\th_t+\s'_t)\si v^{-3}(v^2+v^4)(\c_t(\th_z^\sha)+\c_t(\s'{}_z^\sha)),$$
$$\z v^{-3}(v^6+v^8)(\th_t+\s_t)\si v^{-3}(v^2+v^4)(\c_t(\th_z^\sha)+\c_t(\s_z^\sha)),$$
$$\z v^{-1}(v^4+v^6)(\r_t+\s'_t)\si v\i(1+v^2)(\c_t(\r_z^\sha)+\c_t(\s'{}_z^\sha)),$$
$$\z v^{-1}(v^4+v^6)(\r_t+\s_t)\si v\i(1+v^2)(\c_t(\r_z^\sha)+\c_t(\s_z^\sha)).$$

If $(t,z)$ are as in 4.8 then
$$\z v^{-3}((v^8+v^{10})(\r_t+\s_t)+(v^6+v^8)(\th_t+\s'_t))\si v^{-3}(v^2+v^4)
(\c_t(\th_z^\sha)+\c_t(\s'{}_z^\sha)),$$
$$\z v^{-3}(v^4+v^6)(v^2+v^4)(\th_t+\s_t)\si v^{-3}(v^2+v^4)(\c_t(\th_z^\sha)+\c_t(\s_z^\sha)),$$
$$\z v^{-1}(v^4+v^6)(1+v^2)(\r_t+\s'_t)\si v\i(1+v^2)(\c_t(\r_z^\sha)+\c_t(\s'{}_z^\sha)),$$
$$\z v^{-1}((v^6+v^8)(\th_t+\s'_t)+(v^4+v^6)(\r_t+\s_t))\si v\i(1+v^2)(\c_t(\r_z^\sha)+\c_t(\s_z^\sha)).$$

If $(t,z)$ are as in 4.9 then
$$\z v^{-3}((v^6+v^8)(\r_t+\s_t)+(v^4+v^6)(\th_t+\s'_t))\si v^{-3}(v^2+v^4)
(\c_t(\th_z^\sha)+\c_t(\s'{}_z^\sha)),$$
$$\z v^{-3}(v^4+v^6)(1+v^2)(\r_t+\s'_t)\si v^{-3}(v^2+v^4)(\c_t(\th_z^\sha)+\c_t(\s_z^\sha)),$$
$$\z v^{-1}(1+v^2)(v^2+v^4)(\th_t+\s_t)\si v\i(1+v^2)(\c_t(\r_z^\sha)+\c_t(\s'{}_z^\sha)),$$
$$\z v^{-1}((v^4+v^6)(\th_t+\s'_t)+(v^2+v^4)(\r_t+\s_t))\si v\i(1+v^2)(\c_t(\r_z^\sha)+\c_t(\s_z^\sha)).$$

If $(t,z)$ are as in 4.10 then
$$0\si v^{-3}(v^2+v^4)
(\c_t(\th_z^\sha)+\c_t(\s'{}_z^\sha)),$$
$$0\si v^{-3}(v^2+v^4)(\c_t(\th_z^\sha)+\c_t(\s_z^\sha)),$$
$$0\si v\i(1+v^2)(\c_t(\r_z^\sha)+\c_t(\s'{}_z^\sha)),$$
$$0\si v\i(1+v^2)(\c_t(\r_z^\sha)+\c_t(\s_z^\sha)).$$

\subhead 4.17\endsubhead
Let $z,t\in N_J$. Let ${}^t\fK_J^+=\sum_{C'\in\uCS(L)}\NN[v,v\i]C'_t.$ From 3.5(d) we have (for 
$C\in\uCS(L)$):
$$\c_t(C_z^\sha)=\sum_{C'\in\uCS(L),j'\in\ZZ}\dim{}'V_{t,z,j'}^{C',C}v^{j'}C'_t$$
hence
$$\c_t(C_z^\sha)\in{}^t\fK_J^+.\tag a$$
Using (a) we can extract from the formulas in 4.16 the following facts about $\c_t(C_z^\sha)$. (As in 4.16 
we set $\z=v^{-l(z)+l(t)}$.)

If $(t,z)$ are as in 4.4 then
$$\c_t(\r_z^\sha)\si\z\r_t,\qua \c_t(\s_z^\sha)\si\z\s_t,\qua \c_t(\s'{}_z^\sha)\si\z\s'_t,\qua 
\c_t(\th_z^\sha)\si\z\th_t.$$
If $(t,z)$ are as in 4.5 then
$$\c_t(\r_z^\sha)\si\z v^2\s_t,\qua\c_t(\s_z^\sha)\si\z v^2\r_t,\qua \c_t(s'{}_z^\sha)\si\z v^2\th_t,\qua
\c_t(\th_z^\sha)\si\z v^2\s'_t.$$
If $(t,z)$ are as in 4.6 then   
$$\align&\c_t(\r_z^\sha)\si\z(v^2+v^4)\s_t,\qua \c_t(\s_z^\sha)\si\z(v^2+v^4)\r_t,\qua \\&
\c_t(\s'{}_z^\sha)\si\z (v^2+v^4)\th_t,\qua\c_t(\th_z^\sha)\si\z (v^2+v^4)\s'_t.\endalign$$
If $(t,z)$ are as in 4.7 then
$$\c_t(\r_z^\sha)\si\z v^4\r_t,\qua\c_t(\s_z^\sha)\si\z v^4\s_t,\qua \c_t(\s'{}_z^\sha)\si\z v^4\s'_t, 
\c_t(\th_z^\sha)\si\z v^4\th_t.$$
If $(t,z)$ are as in 4.8 then
$$\align&\c_t(\r_z^\sha)\si\z(v^4\r_t+v^6\s'_t),\qua\c_t(\s_z^\sha)\si\z(v^4\s_t+v^6\th_t),\qua\\&
\c_t(\s'{}_z^\sha)\si\z(v^4\s'_t+v^6\r_t),\qua\c_t(\th_z^\sha)\si\z(v^4\th_t+v^6\s_t).\endalign$$
If $(t,z)$ are as in 4.9 then
$$\align&\c_t(\r_z^\sha)\si\z(v^2\s_t+v^4\th_t),\qua\c_t(\s_z^\sha)\si\z(v^2\r_t+v^4\s'_t),\qua\\&
\c_t(\s'{}_z^\sha)\si\z(v^2\th_t+v^4\s_t),\qua \c_t(\th_z^\sha)\si\z(v^2\s'_t+v^4\r_t).\endalign$$
If $(t,z)$ are as in 4.10 then
$$\c_t(\r_z^\sha)\si0,\qua\c_t(\s_z^\sha)\si0,\qua\c_t(\s'{}_z^\sha)\si0,\qua \c_t(\th_z^\sha)\si0.$$

\subhead 4.18\endsubhead
Let $z,t\in N_J$. Using the results in 4.17 we see that   
$$\c_t(\r_z^\sha)-\c_t(\s_z^\sha)-\c_t(\s'{}_z^\sha)+\c_t(\th_z^\sha)\si X_{t,z}(\r_t-\s_t-\s'_t+\th_t),
\tag a$$
where $X_{t,z}\in\ca$ is as follows:

$X_{t,z}=v^{-l(z)+l(t)} $ if $(t,z)$ are as in 4.4;

$X_{t,z}=-v^{-l(z)+l(t)} v^2$ if $(t,z)$ are as in 4.5;

$X_{t,z}=-v^{-l(z)+l(t)} (v^2+v^4)$ if $(t,z)$ are as in 4.6;

$X_{t,z}=v^{-l(z)+l(t)}  v^4$ if $(t,z)$ are as in 4.7;

$X_{t,z}=v^{-l(z)+l(t)} (v^4-v^6)$ if $(t,z)$ are as in 4.8;

$X_{t,z}=-v^{-l(z)+l(t)}(v^2-v^4)$ if $(t,z)$ are as in 4.9;

$X_{t,z}=0$ if $(t,z)$ are as in 4.10.

It follows that

$X_{t,z}=\e(z)\e(t)v^{-\cl(z)+\cl(t)} $ if $(t,z)$ are as in 4.4, 4.5 or 4.7;

$X_{t,z}=\e(z)\e(t)v^{-\cl(z)+\cl(t)} (1+v^2)$ if $(t,z)$ are as in 4.6;

$X_{t,z}=\e(z)\e(t)v^{-\cl(z)+\cl(t)} (1-v^2)$ if $(t,z)$ are as in 4.8 or 4.9;

$X_{t,z}=0$ if $(t,z)$ are as in 4.10.

\subhead 4.19\endsubhead
Define $\ch,p_{t,z}$ as in 0.1 in terms of $\cw=N_J$, $\cl$. According to \cite{\HEC, 7.6} we have:

(i) $p_{t,z}=v^{-\cl(z)+\cl(t)} (1+v^2)$ if $z=efe$ and $t\in\{\emp,e\}$;

(ii) $p_{t,z}=v^{-\cl(z)+\cl(t)} (1-v^2)$ if $z=fef$ and $t\in\{\emp,f\}$;

(iii) $p_{t,z}=v^{-\cl(z)+\cl(t)}$ if $t\le z$ in the usual partial order of $N_J$ with $(t,z)$ not as
in (i),(ii);

(iv) $p_{t,z}=0$ if $t\not\le z$.
\nl
We can now restate the result in 4.18 as follows.
$$\c_t(\r_z^\sha)-\c_t(\s_z^\sha)-\c_t(\s'{}_z^\sha)+\c_t(\th_z^\sha)\si 
p_{t,z}\e(z)\e(t)(\r_t-\s_t-\s'_t+\th_t).$$
We see that Conjectures 1,2,3 in 3.11 hold in our case.


\widestnumber \key{BBD}
\Refs
\ref\key\BE\by R.B\'edard\paper On the Brauer liftings for modular representations\jour J.Algebra\vol93\yr
1985\pages332-353\endref
\ref\key\BBD\by A.Beilinson,J.Bernstein and P.Deligne\book Faisceaux pervers\bookinfo Analysis and topology 
on singular spaces, I\jour Ast\'erisque\vol100\publ Soc. Math. France, Paris\yr1982\pages5-171\endref
\ref\key\BRA\by T.Braden\paper Hyperbolic localization of intersection cohomology\jour Transfor.Groups\vol8
\yr2003\pages209-216\endref
\ref\key\GOR\by M.Goresky\paper Tables of Kazhdan-Lusztig polynomials\jour 
www.math.ias.edu/~goresky\yr1996\endref
\ref\key\HE\by X.He\paper The $G$-stable pieces of the wonderful compactification\jour Trans.Amer.Math.Soc.
\vol359\yr2007\pages3005-3024\endref
\ref\key\KL\by D.Kazhdan and G.Lusztig\paper Representations of Coxeter groups and Hecke algebras\jour
Invent. Math.\vol53\yr1979\pages165-184\endref
\ref\key\KLL\by D.Kazhdan and G.Lusztig\paper Schubert varieties and Poincar\'e duality\jour
Proc. Symp. Pure Math.\vol36\paperinfo Amer. Math. Soc.\yr1980\pages185-203\endref 
\ref\key\COX\by G.Lusztig\paper Coxeter orbits and eigenspaces of Frobenius
\jour Inv.Math.\vol28\yr1976\pages101-159\endref
\ref\key\CLA\bysame\paper Irreducible representations of finite classical groups
\jour Inv.Math.\vol43\yr1977\pages125-175\endref
\ref\key\ORA\bysame\book Characters of reductive groups over a finite field\bookinfo
 Ann.Math.Studies 107\publ Princeton U.Press\yr1984\endref
\ref\key\CLASSP\bysame\paper Classification of unipotent representations of simple $p$-adic groups
\jour Int. Math. Res. Notices\yr1995\pages517-589\endref
\ref\key\HEC\bysame\book Hecke algebras with unequal parameters\bookinfo CRM Monograph Ser. 18, 
Amer. Math. Soc.\yr2003\endref
\ref\key\PCS\bysame\paper Parabolic character sheaves,I\jour Moscow Math.J.\vol4\yr2004\pages153-179
\endref
\ref\key\CDGVIII\bysame\paper Character sheaves on disconnected groups,VIII\jour Represent.Theory
\vol10\yr2006\pages314-352\endref
\ref\key\PCSIII\bysame\paper Parabolic character sheaves,III\jour Moscow Math.J.\vol10\yr2010\pages603-609
\endref
\ref\key\URCC\bysame\paper Unipotent representations as a categorical center\jour arxiv:1401.2889
\endref
\endRefs
\enddocument